\documentclass{article}

\usepackage{graphicx,verbatim,array,multicol,courier}
\usepackage{amsmath,amsthm,amssymb}
\usepackage[utf8]{inputenc}
\usepackage{color, multicol}

\usepackage{enumitem}

\usepackage{subcaption}

\usepackage{natbib}

\usepackage{bm}
\usepackage{bbm}

\usepackage{url}

\usepackage[left=4cm,right=4cm,top=3cm,bottom=3cm]{geometry}

\DeclareFontFamily{U}{mathx}{\hyphenchar\font45}
\DeclareFontShape{U}{mathx}{m}{n}{
      <5> <6> <7> <8> <9> <10>
      <10.95> <12> <14.4> <17.28> <20.74> <24.88>
      mathx10
      }{}
\DeclareSymbolFont{mathx}{U}{mathx}{m}{n}
\DeclareFontSubstitution{U}{mathx}{m}{n}
\DeclareMathAccent{\widecheck}{0}{mathx}{"71}

\definecolor{navy}{rgb}{0,0,0.502}
\definecolor{brown}{rgb}{0.59, 0.29, 0.0}
\definecolor{green}{rgb}{0.25,0.50,0.25}
\definecolor{violet}{rgb}{0.56, 0.0, 1.0}

\newcommand{\ind}{\mathbbm{1}}

\newtheorem{Theo}{Theorem}

\newtheorem{Coro}{Corollary}

\theoremstyle{definition}

\newtheorem{Rmk}{Remark} 

\title{Optimal pooling and distributed inference for the tail index and extreme quantiles}
\author{Abdelaati Daouia$^{a}$, Simone A. Padoan$^{b}$ \& Gilles Stupfler$^{c}$}
\date{$^{a}$ {\small Toulouse School of Economics, University of Toulouse Capitole, France} \\ $^{b}$ {\small Department of Decision Sciences, Bocconi University of Milan, via Roentgen 1, 20136 Milano, Italy} \\ $^{c}$ {\small Univ Rennes, Ensai, CNRS, CREST - UMR 9194, F-35000 Rennes, France}}

\begin{document}

\maketitle

\begin{abstract}
This paper investigates pooling strategies for tail index and extreme quantile estimation from heavy-tailed data. To fully exploit the information contained in several samples, we present general weighted pooled Hill estimators of the tail index and weighted pooled Weissman estimators of extreme quantiles calculated through a nonstandard geometric averaging scheme. We develop their large-sample asymptotic theory across a fixed number of samples, covering the general framework of heterogeneous sample sizes with different and asymptotically dependent distributions. Our results include optimal choices of pooling weights based on asymptotic variance and MSE minimization. In the important application of distributed inference, we prove that the variance-optimal distributed estimators are asymptotically equivalent to the benchmark Hill and Weissman estimators based on the unfeasible combination of subsamples, while the AMSE-optimal distributed estimators enjoy a smaller AMSE than the benchmarks in the case of large bias. We consider additional scenarios where the number of subsamples grows with the total sample size and effective subsample sizes can be low. We extend our methodology to handle serial dependence 
and the presence of covariates. Simulations confirm that our pooled estimators perform virtually as well as the benchmark estimators. Two applications to real weather and insurance data are showcased.
\end{abstract}



\section{Introduction} 

The contemporary problem of efficient analysis of massive data has led to the development of 
the divide-and-conquer approach,
which consists in dividing data into multiple samples that can be processed across several machines, before combining the results from subsamples on a central machine 
by making use of pooling techniques. In statistics, this strategy gives rise to distributed inference from many datasets, allowing to alleviate the computational challenges and constraints in storage imposed by the availability of extremely large datasets, or to handle data privacy issues as in banking and insurance. Pioneering contributions in the distributed framework focused on the regression mean and central parameters. In the last three years, distributed inferential procedures have been also developed for quantile estimation, mainly in a regression setup. Prominent among these contributions are~\cite{volchache2019},~\cite{xucaijiasunhua2020}, and~\cite{wanma2021}. The analyses therein concentrate on the statistical inference of ordinary conditional quantiles at fixed tail probability levels.

The estimation of extreme quantiles in a distributed computing setting, when their order tends to 1 as the total sample size goes to infinity, has not yet received any attention, however. More generally, distributed inference for tail quantities from the perspective of extreme value theory is still in its infancy. To the best of our knowledge, only tail index estimation has been recently explored in~\cite{chelizho2021} by taking a simple average of subsample tail index estimators as the final distributed estimator. 
From a broader perspective, the use of machine learning methods for extreme value analysis has started only very recently with, {\it e.g.},~\cite{ahmeinzho2021} for extreme value statistics in semi-supervised models,~\cite{veldomcaieng2021} for a gradient boosting procedure to estimate extreme conditional quantiles, and~\cite{aghporsabzho2021} for dimension reduction based on tail inverse regression.

In this article we go further than~\cite{chelizho2021} by addressing several important questions about distributed inference for the tail index, but also for extreme quantiles of heavy-tailed data. Our approach is based on the use of a general terminology and theory of pooling that encompasses distributed estimation from different samples whose distributions have a common target parameter, given in our setup by the tail index or an extreme quantile. Instead of naively averaging the 
estimators calculated 
from the available subsamples, it is of interest, both from a theoretical and a practical perspective, to construct a general class of weighted pooled estimators and to establish a fully data-driven inferential procedure integrating the optimal choice of weights. In particular, our overarching goal is to develop the asymptotic theory of the optimally pooled tail index and extreme quantile 
estimators, under weak technical conditions, covering both scenarios where the number of available subsamples 
is bounded or growing with the total sample size, as well as the general situation where the observed data can be dependent within and/or across subsamples. The pooling approach itself has a rich history dating back to~\cite{coc1937} for the estimation of the common mean of several samples. It is also encountered in the econometrics of panel data since the influential work of~\cite{the1954} and~\cite{zel1962}. The idea of combining pooling and extreme value techniques has originally been suggested by~\cite{kinfrilil2016},~\cite{asaengdav2018}, and~\cite{vigengzsc2021} in climate science problems.

The first contribution of this paper is a joint asymptotic normality result for Hill estimators~\citep{hil1975} calculated from a fixed number $m$ of samples with heavy-tailed data. In particular, we allow for different sample sizes, effective sample sizes, marginal distributions, 
and for dependence across samples.
We apply this general result to design optimal pooling strategies of subsample Hill estimators for tail index estimation. We consider optimal weights that minimize either the asymptotic variance or the Asymptotic Mean Squared Error (AMSE) of the pooled estimator. These developments rely on a very general theory that we derive for a generic weighted pooled estimator built from $m$ subsample estimators for a common unknown parameter. This theory comes into play when the subsample estimators are jointly asymptotically normal, and can be biased and correlated. To the best of our knowledge, no such unrestricted approach has been fully investigated. We also construct bias-reduced versions of the proposed pooled tail index estimators. Then we discuss the problem of ultimate interest in extreme value analysis, which lies in estimating extreme quantiles either locally in each machine in the tail homogeneous setting of equal marginal tail indices, where tail quantiles are possibly only asymptotically proportional across subsamples, or globally by pooling subsample extreme quantile Weissman estimators~\citep{wei1978} in the more restrictive tail homoskedastic setting of asymptotically equivalent marginal tail quantiles. Our approach in both cases relies on a specific weighted geometric pooling scheme, particularly relevant for extreme quantiles, as opposed to arithmetic averaging naturally used for pooling ordinary quantiles~\citep{knibas2003}. Moreover, we explore inferential aspects of pooling for extreme values by constructing likelihood ratio-type tests for either tail homogeneity or tail homoskedasticity, as well as asymptotic confidence intervals for the tail index and extreme quantiles.

We also specialize the discussion to distributed inference as an important application of our general 
theory. In this particular case, 
due to computational costs or privacy restrictions, 
the data 
in sample number $j$ can only be processed by the $j$th machine, with very restricted or no communication allowed between the $m$ machines, before an end user operating from a central machine conducts distributed inference from limited information transmitted by each machine. It is also assumed that the data is independent and identically distributed (i.i.d.)~within and across machines. Under this setup, we examine and compare the asymptotic theory of the distributed tail index and extreme quantile estimators to the behavior of their respective benchmark Hill and Weissman estimators based on the unfeasible direct combination of subsamples. We extend this theory further by considering first the case when effective sample sizes are highly unbalanced among machines, and then the case of a growing number of machines~$m=m(n)$ with the total sample size~$n$. Finally, we tackle the problem of serial dependence within the data in the presence of covariates, showing how appropriately filtering the observations allows to recover the asymptotic theory from independent observations.

What first distinguishes our contribution relative to earlier literature, and in particular to~\cite{chelizho2021} for tail index estimation, is that a detailed study is conducted for the case of a fixed number $m$ of subsamples with different sample sizes. These considerations are motivated by practical concerns. For example, the financial application in~\cite{chelizho2021} itself requires the rather small value $m=5$ and a ratio between lowest and highest subsample sizes roughly equal to $25\%$, whereas their theory only considers the case $m=m(n)\to\infty$ as $n\to\infty$, with equal subsample sizes. 
We allow for general choices of weights in the pooling scheme, we construct bias-reduced versions of our pooled estimators, and we design inference procedures under very weak conditions that hold for any reasonable heavy-tailed model. 
We construct asymptotic variance-optimal and AMSE-optimal estimators, 
with meaningful comparisons between these optimally pooled estimators and the unfeasible benchmark Hill estimator. We 
revisit the case $m=m(n)\to\infty$ by providing, unlike~\cite{chelizho2021}, a unified convergence result for tail index estimation that handles both bounded and unbounded effective sample sizes in the marginal Hill estimators. In particular, we carefully derive the correct expression of the asymptotic variance for the distributed Hill estimator, which should be different from the expression provided by~\cite{chelizho2021} in the case of unbalanced effective sample sizes. 
Moreover, this is the first work to implement the idea of geometric weighted extreme quantile pooling; our experience with simulated and real data indicates the superiority of the corresponding estimators over the arithmetically pooled competitors that suffer from substantially larger bias and variance.

The paper is organized as follows. Section~\ref{sec:main} develops our general pooling theory for tail index and extreme quantile estimation, while Section~\ref{sec:distrinf} focuses on the special framework of distributed inference. Section~\ref{sec:filtering} extends our methodology to handle serial dependence and the presence of covariates through filtering. Section~\ref{sec:fin} illustrates the usefulness of the proposed extreme value pooling and distributed inference methods through a simulation study and concrete applications to weather and insurance data. The supplement to this article contains additional theoretical results and all the proofs, with further details on our simulation study. Our methods and data have been incorporated into the open-source {\tt R} package {\tt ExtremeRisks}.

\section{Pooling extreme value estimators}
\label{sec:main}
\subsection{Pooled Hill estimators of the tail index}
\label{sec:main:EVI}

Let $\boldsymbol{X}=(X_1,\ldots,X_m)^{\top}$ denote an $m-$dimensional random vector, and $\boldsymbol{X}_i=(X_{i,1},\ldots,X_{i,m})^{\top}$ ($i\geq 1$) denote independent copies of $\boldsymbol{X}$. We assume that the available data consists of the $X_{i,j}$, for $1\leq j\leq m$ and $1\leq i\leq n_j=n_j(n)$, with $n=\sum_{j=1}^m n_j$ being the total number of univariate data points available across all samples. We suppose that $n_j\to\infty$ as $n\to\infty$. Table~\ref{tab:data} provides a simple representation of the available data when $n_1<n_2<\hdots<n_m$.
\begin{table}[htbp]
\centering
\begin{tabular}{c|ccccccccc}
Vector $\boldsymbol{X}$ & \multicolumn{9}{c}{Available data}\\
\hline \\
$X_1$ & $X_{1,1}$ & $\hdots$ & $X_{n_1,1}$ & & & & & &  \\
$X_2$ & $X_{1,2}$ & $\hdots$ & $X_{n_1,2}$ & $X_{n_1+1,2}$ & $\hdots$ & $X_{n_2,2}$ & & & \\
$\vdots$ &  \\
$X_m$ & $X_{1,m}$ & $\hdots$ & $X_{n_1,m}$ & $X_{n_1+1,m}$ & $\hdots$ & $X_{n_2,m}$ & $X_{n_2+1,m}$ & $\hdots$ & $X_{n_m,m}$ \\
\end{tabular}
\caption{\label{tab:data} Representation of the available data when $n_1<n_2<\hdots<n_m$. The random vectors $\boldsymbol{X}_i=(X_{i,1},\ldots,X_{i,m})^{\top}$ ($i\geq 1$) are independent copies of the vector $\boldsymbol{X}=(X_1,\ldots,X_m)^{\top}$.}
\end{table}

We focus on the general framework where the components $X_j$ of the random vector $\bm{X}$ have continuous, right heavy-tailed distribution functions $F_j$, with associated survival functions $\overline{F}_j=1-F_j$ and tail quantile functions $U_j:t\mapsto \inf\{ x\in \mathbb{R} \, | \, 1/\overline{F}_j(x) \geq t \}$ that satisfy
\vskip1ex
\noindent
$\mathcal{C}_2(\boldsymbol{\gamma},\boldsymbol{\rho},\boldsymbol{A})$ For any $j\in \{1,\ldots,m\}$, the function $U_j$ satisfies the second-order condition:
\vskip1ex
\noindent
$\mathcal{C}_2(\gamma_j,\rho_j,A_j)$ 
$U_j$ is second-order regularly varying in a neighborhood of $+\infty$ with index $\gamma_j>0$, second-order parameter $\rho_j\leq 0$ and an auxiliary function $A_j$ having constant sign and converging to 0 at infinity, that is,
$$
\forall x>0, \ \lim_{t\to \infty}\frac{1}{A_j(t)} \left[ \frac{U_j(tx)}{U_j(t)} - x^{\gamma_j} \right] = x^{\gamma_j} \frac{x^{\rho_j}-1}{\rho_j}
$$
where the right-hand side should be read as $x^{\gamma_j} \log x$ when $\rho_j=0$. 
\vskip1ex
\noindent
In this condition $|A_j|$ is regularly varying with index $\rho_j$~\citep[by Theorems~2.3.3 and~2.3.9 in][]{haafer2006}, meaning that the larger $|\rho_j|$ is, 
the smaller the error in the approximation of the right tail of $U_j$ by a purely Pareto tail will be. All usual heavy-tailed distributions 
satisfy these conditions, see Table 2.1 in~p.59 of~\cite{beigoesegteu2004} for a detailed list of examples.

To incorporate the dependence between samples into the inference procedure, we assume an appropriate pairwise tail dependence structure based on the functions
$\overline{C}_{j,\ell}(u,v) = \mathbb{P}( \overline{F}_j(X_j)\leq u, \ \overline{F}_{\ell}(X_\ell)\leq v )$, for $u,v\in [0,1]$ that are essentially the bivariate survival copulae of 
$\bm{X}$, namely:
\vskip1ex
\noindent
$\mathcal{J}(\boldsymbol{R})$ For any $(j,\ell)$ with $j\neq \ell$, there is a function $R_{j,\ell}$ on $[0,\infty]^2 \setminus \{  (\infty,\infty)\}$ such that
\[
\forall (x_j,x_{\ell})\in [0,\infty]^2 \setminus \{  (\infty,\infty)\}, \ \lim_{s\to \infty} s \, \overline{C}_{j,\ell}(x_j/s,x_\ell/s) = R_{j,\ell}(x_j,x_\ell). 
\]
This condition imposes the existence of a limiting dependence structure in the joint right tail of $X_j$ and $X_{\ell}$, given by the {\it tail copula} $R_{j,\ell}$~(see~\cite{schsta2006}). It can be viewed as a minimal assumption when it comes to assessing the dependence structure between extreme value estimators. 

After ordering the data in the $j$th sample as $X_{1:n_j,j}\leq X_{2:n_j,j}\leq \cdots \leq X_{n_j:n_j,j}$, (the notation is inspired by~\cite{davnag2003}), 
we introduce the marginal Hill estimators 
$
\widehat{\gamma}_j(k_j) = k_j^{-1} \sum_{i=1}^{k_j} \log(X_{n_j-i+1:n_j,j}/X_{n_j-k_j:n_j,j})
$
which involve the top $(k_j + 1)$ highest order statistics in each sample, for $k_j=k_j(n)\geq 1$. 
The integer $k_j$ is the {\it effective} sample size in sample $j$, and we set $k=\sum_{j=1}^m k_j$ to be the {\it total effective} sample size in the vector of estimators $\widehat{\boldsymbol{\gamma}}_n=\widehat{\boldsymbol{\gamma}}_n(k_1,\ldots,k_m)=(\widehat{\gamma}_1(k_1),\ldots,\widehat{\gamma}_m(k_m))^{\top}$. 
Our ultimate interest is in the case $\boldsymbol{\gamma}=\gamma\boldsymbol{1}$ where the $\gamma_j$ are equal to a common 
$\gamma$ estimated by
\[
\widehat{\gamma}_n(\boldsymbol{\omega}) = \widehat{\gamma}_n(\omega_1,\ldots,\omega_m) = \sum_{j=1}^m \omega_j \widehat{\gamma}_j(k_j) = \boldsymbol{\omega}^{\top} \widehat{\boldsymbol{\gamma}}_n, \mbox{ with } \boldsymbol{\omega}^{\top} \boldsymbol{1} = 1.
\]
The asymptotic distribution of any element within this class of estimators is stated in the following theorem, along with the joint asymptotic normality of $\widehat{\gamma}_j(k_j)$, for $j=1,\ldots,m$.
\begin{Theo}
\label{theo:pool_mfixed_kinfty_general}
Assume that conditions $\mathcal{C}_2(\boldsymbol{\gamma},\boldsymbol{\rho},\boldsymbol{A})$ and $\mathcal{J}(\boldsymbol{R})$ hold. Suppose, without loss of generality, that $n_1\leq n_2\leq\cdots\leq n_m$, and then that $k_j=k_j(n)\to\infty$ with $k_j/n_j\to 0$, $n_1/n_j\to b_j\in (0,1]$, $k_1/k_j\to c_j\in (0,\infty)$ (with $c_1=1$ and $b_m\leq b_{m-1}\leq\cdots\leq b_1=1$) and $\sqrt{k_j} A_j(n_j/k_j)\to \lambda_j\in\mathbb{R}$ for any $j\in \{1,\ldots,m\}$, as $n\to\infty$. Let the weight vector $\boldsymbol{\omega} = (\omega_1,\ldots,\omega_m)^{\top}$ be such that $\boldsymbol{\omega}^{\top} \boldsymbol{1}=1$ and define a 
vector $\boldsymbol{B}$ and 
symmetric matrix $\mathbf{V}$ by
\[
\boldsymbol{B} = \left( \frac{\lambda_1}{1-\rho_1},\ldots, \frac{\lambda_m}{1-\rho_m} \right)^{\top} \mbox{ and } \mathbf{V}_{j,\ell} = \begin{cases} \gamma_j^2  & \mbox{if } j=\ell, \\[5pt] \gamma_j \gamma_{\ell} \sqrt{\dfrac{c_j}{c_{\ell}}} R_{j,\ell}(c_{\ell}/c_j,b_{\ell}/b_j) & \mbox{if } j<\ell. \end{cases} 
\]
Then 
$
(\sqrt{k_1} (\widehat{\gamma}_1(k_1)-\gamma_1),\ldots,\sqrt{k_m}(\widehat{\gamma}_m(k_m)-\gamma_m))^{\top} \stackrel{d}{\longrightarrow} \mathcal{N}(\boldsymbol{B},\mathbf{V}).
$
In particular, if $\boldsymbol{\gamma}=\gamma\boldsymbol{1}$, then
\begin{align*}
    & \sqrt{k} (\widehat{\gamma}_n(\boldsymbol{\omega})-\gamma) \stackrel{d}{\longrightarrow} \mathcal{N}\left( \boldsymbol{\omega}^{\top} \boldsymbol{B}_{\boldsymbol{c}}, \boldsymbol{\omega}^{\top} \mathbf{V}_{\boldsymbol{c}} \boldsymbol{\omega} \right) , \\
    \mbox{with } & \boldsymbol{B}_{\boldsymbol{c}} = \sqrt{ \sum_{i=1}^m c_i^{-1} } \left( \sqrt{c_1}\frac{\lambda_1}{1-\rho_1},\ldots, \sqrt{c_m}\frac{\lambda_m}{1-\rho_m} \right)^{\top} \\
    \mbox{and } & [\mathbf{V}_{\boldsymbol{c}}]_{j,\ell} = \left(\sum_{i=1}^m c_i^{-1}\right) \begin{cases} \gamma^2 c_j  & \mbox{if } j=\ell, \\[5pt] \gamma^2 c_j R_{j,\ell}( c_{\ell}/c_j,b_{\ell}/b_j ) & \mbox{if } j<\ell. \end{cases} 
\end{align*}
The matrix $\mathbf{V}$ is positive definite if and only if $\mathbf{V}_{\boldsymbol{c}}$ is so, and hence we have the following results on optimal weights: 
\begin{enumerate}[wide, labelindent=0pt]
\item (Variance-optimal weights) There is a unique solution to the minimization problem of $\boldsymbol{\omega}^{\top}\mathbf{V}_{\boldsymbol{c}}\boldsymbol{\omega}$ subject to the constraint $\boldsymbol{\omega}^{\top} \boldsymbol{1} = 1$, which is 
\[
\boldsymbol{\omega}^{(\mathrm{Var})} = \frac{\mathbf{V}_{\boldsymbol{c}}^{-1} \boldsymbol{1}}{\boldsymbol{1}^{\top} \mathbf{V}_{\boldsymbol{c}}^{-1} \boldsymbol{1}}, \mbox{ and then } \sqrt{k} \left(\widehat{\gamma}_n(\boldsymbol{\omega}^{(\mathrm{Var})})-\gamma\right) \stackrel{d}{\longrightarrow} \mathcal{N}\left( \frac{\boldsymbol{1}^{\top} \mathbf{V}_{\boldsymbol{c}}^{-1} \boldsymbol{B}_{\boldsymbol{c}}}{\boldsymbol{1}^{\top} \mathbf{V}_{\boldsymbol{c}}^{-1} \boldsymbol{1}}, \frac{1}{\boldsymbol{1}^{\top} \mathbf{V}_{\boldsymbol{c}}^{-1} \boldsymbol{1}} \right).
\]
\item (AMSE-optimal weights) There is a unique solution to the minimization problem of $\mathrm{AMSE}(\boldsymbol{\omega}) = (\boldsymbol{\omega}^{\top}\boldsymbol{B}_{\boldsymbol{c}})^2+\boldsymbol{\omega}^{\top}\mathbf{V}_{\boldsymbol{c}}\boldsymbol{\omega}$ subject to the constraint $\boldsymbol{\omega}^{\top} \boldsymbol{1} = 1$, which is 
\[
\boldsymbol{\omega}^{(\mathrm{AMSE})} = \frac{(1+\boldsymbol{B}_{\boldsymbol{c}}^{\top} \mathbf{V}_{\boldsymbol{c}}^{-1}\boldsymbol{B}_{\boldsymbol{c}}) \mathbf{V}_{\boldsymbol{c}}^{-1} \boldsymbol{1} - (\boldsymbol{1}^{\top} \mathbf{V}_{\boldsymbol{c}}^{-1} \boldsymbol{B}_{\boldsymbol{c}}) \mathbf{V}_{\boldsymbol{c}}^{-1}\boldsymbol{B}_{\boldsymbol{c}}}{(1+\boldsymbol{B}_{\boldsymbol{c}}^{\top} \mathbf{V}_{\boldsymbol{c}}^{-1}\boldsymbol{B}_{\boldsymbol{c}}) (\boldsymbol{1}^{\top} \mathbf{V}_{\boldsymbol{c}}^{-1} \boldsymbol{1}) - (\boldsymbol{1}^{\top} \mathbf{V}_{\boldsymbol{c}}^{-1} \boldsymbol{B}_{\boldsymbol{c}})^2}.
\]
The optimal value of $\mathrm{AMSE}(\boldsymbol{\omega})$ is 
\[
\mathrm{AMSE}(\boldsymbol{\omega}^{(\mathrm{AMSE})}) = \frac{1+\boldsymbol{B}_{\boldsymbol{c}}^{\top} \mathbf{V}_{\boldsymbol{c}}^{-1}\boldsymbol{B}_{\boldsymbol{c}}}{(1+\boldsymbol{B}_{\boldsymbol{c}}^{\top} \mathbf{V}_{\boldsymbol{c}}^{-1}\boldsymbol{B}_{\boldsymbol{c}}) (\boldsymbol{1}^{\top} \mathbf{V}_{\boldsymbol{c}}^{-1} \boldsymbol{1}) - (\boldsymbol{1}^{\top} \mathbf{V}_{\boldsymbol{c}}^{-1} \boldsymbol{B}_{\boldsymbol{c}})^2}.
\]
\end{enumerate}
Finally, if
$\widehat{\boldsymbol{\omega}}_n^{\top} \boldsymbol{1}=1$ with $\widehat{\boldsymbol{\omega}}_n\stackrel{\mathbb{P}}{\longrightarrow} \boldsymbol{\omega}$, 
then the composite estimator $\widehat{\gamma}_n(\widehat{\boldsymbol{\omega}}_n)$ is $\sqrt{k}-$asymptotically equivalent to $\widehat{\gamma}_n(\boldsymbol{\omega})$ in the sense that $\sqrt{k}(\widehat{\gamma}_n(\widehat{\boldsymbol{\omega}}_n) - \widehat{\gamma}_n(\boldsymbol{\omega}))=\operatorname{o}_{\mathbb{P}}(1)$.
\end{Theo}
The proof requires 
applying a very general pooling result in conjunction with the assumption that the sample sizes $n_j$ are asymptotically proportional (possibly unbalanced) and so are the effective sample sizes $k_j$. This ensures that none of the $\widehat{\gamma}_j(k_j)$ imposes its limiting distribution to the others. 
\begin{Rmk}[Bias components across samples] 
\label{rmk:bias_sites}
The structure of the asymptotic bias component $\boldsymbol{B}_{\boldsymbol{c}}$ is constrained by the second-order parameters $\rho_j$: for all its components to be non-zero, the $\rho_j$ should be equal, because of the proportionality assumption on the $k_j$'s and the $n_j$'s with 
the fact that the functions $|A_j|$ are regularly varying with index $\rho_j$. More specifically,
if $\rho^* = \max_{1\leq j \leq m} \rho_j$, then the $j$th component of $\boldsymbol{B}_{\boldsymbol{c}}$ is 0 whenever $\rho_j<\rho^*$. In the important special context of distributed inference (see Section~\ref{sec:distrinf}), where all marginal distributions are the same, an asymptotic bias component will be present in every marginal.
\end{Rmk}
\subsection{Optimal choices of weights}
\label{sec:main:optimal_EVI}
The idea now is to use the pooled estimator $\widehat{\gamma}_n(\boldsymbol{\omega})$ in conjunction with the optimal choices of the weight vector $\boldsymbol{\omega}$. Since the resulting optimal values of weights depend on the
asymptotic bias and variance components in view of 
Theorem~\ref{theo:pool_mfixed_kinfty_general}, they should be estimated first. One convenient way of estimating the bias component, which hinges on $\lambda_j$ and $\rho_j$, is by assuming $A_j(t)=\gamma \beta_j t^{\rho_j}$ for some constant $\beta_j$. Most commonly used heavy-tailed models satisfy this mild proportionality condition between $A_j(t)$ and $t^{\rho_j}$, see Table~1 in~\cite{girstuauc2021}. Under this assumption, consistent estimators $\widehat{\beta}_j$ and $\widehat{\rho}_j$ of $\beta_j$ and $\rho_j$ are available and implemented in open-source software, for example in the R function {\tt mop} from the package {\tt evt0}. This yields an estimator of 
$\lambda_j = \lim_{n\to\infty} \sqrt{k_j} A_j(n_j/k_j)$ as 
$
\widehat{\lambda}_j = \sqrt{k_j} \times \widehat{\gamma}_n(\boldsymbol{\omega}) \widehat{\beta}_j (n_j/k_j)^{\widehat{\rho}_j}.
$
Here the choice of 
$\boldsymbol{\omega}$ in the $\gamma$ estimator is arbitrary; for example, without any prior knowledge of the dependence structure, the naive weights $\boldsymbol{\omega}=(1/m,\ldots,1/m)^{\top}$ seem reasonable.

Now, to estimate the covariance matrix, let $n_{j,\ell}=\min(n_j,n_{\ell})$ and $k_{j,\ell}=k_j$ if $n_j<n_{\ell}$ and $k_{\ell}$ otherwise, and consider the estimator of the tail copula function $R_{j,\ell}$ defined as
\[
\widehat{R}_{j,\ell}(u,v) = \widehat{R}_{j,\ell}(u,v;k_{j,\ell}) = \frac{1}{k_{j,\ell}}\sum_{i=1}^{n_{j,\ell}}\ind
\left\{
\frac{n_{j,\ell}+1 - r_{n_{j,\ell},i,j}}{k_{j,\ell}(n_{j,\ell}+1)/n_{j,\ell}}\leq u, \frac{n_{j,\ell}+1 - r_{n_{j,\ell},i,\ell}}{k_{j,\ell}(n_{j,\ell}+1)/n_{j,\ell}}\leq v
\right\}.
\]
[Here $r_{n_{j,\ell},i,j}$ (resp.~$r_{n_{j,\ell},i,\ell}$) 
stands for 
the rank of 
$X_{i,j}$ (resp.~$X_{i,\ell}$) among the observations $X_{1,j},X_{2,j},\ldots,X_{n_{j,\ell},j}$ (resp.~$X_{1,\ell},X_{2,\ell},\ldots,X_{n_{j,\ell},\ell}$), namely, the first $n_{j,\ell}$ observations in sample $j$ (resp.~$\ell$).] This is a modified version of the empirical upper tail copula in Equation~(13) of~\cite{schsta2006}. Adapting Lemma~7 from~\cite{stu2019} shows that it is a locally uniformly consistent estimator of $R_{j,\ell}$ on $(0,\infty)^2$ under our technical conditions. 
Combining these tools, we arrive at the estimators 
\begin{align*}
\widehat{\boldsymbol{B}}_{\boldsymbol{c}} &= 
\sqrt{k} \left( \frac{\widehat{\lambda}_1/\sqrt{k_1}}{1-\widehat{\rho}_1}, \frac{\widehat{\lambda}_2/\sqrt{k_2}}{1-\widehat{\rho}_2}, \ldots, \frac{\widehat{\lambda}_m/\sqrt{k_m}}{1-\widehat{\rho}_m} \right)^{\top} \\
\mbox{and } & [\widehat{\mathbf{V}}_{\boldsymbol{c}}]_{j,\ell} = 
k \, \widehat{\gamma}_n^2(\boldsymbol{\omega}) \begin{cases} \dfrac{1}{k_j} & \mbox{if } j=\ell, \\[10pt] \dfrac{1}{k_j} \widehat{R}_{j,\ell}(k_j/k_{\ell},n_j/n_{\ell}) & \mbox{if } j\neq\ell. \end{cases}
\end{align*}
According to 
Theorem~\ref{theo:pool_mfixed_kinfty_general}, the vectors of variance-optimal weights and AMSE-optimal weights can then be estimated by
\[
\widehat{\boldsymbol{\omega}}_n^{(\mathrm{Var})} = \frac{\widehat{\mathbf{V}}_{\boldsymbol{c}}^{-1} \boldsymbol{1}}{\boldsymbol{1}^{\top} \widehat{\mathbf{V}}_{\boldsymbol{c}}^{-1} \boldsymbol{1}} \mbox{ and }  
\widehat{\boldsymbol{\omega}}_n^{(\mathrm{AMSE})} = \frac{(1+\widehat{\boldsymbol{B}}_{\boldsymbol{c}}^{\top} \widehat{\mathbf{V}}_{\boldsymbol{c}}^{-1}\widehat{\boldsymbol{B}}_{\boldsymbol{c}}) \widehat{\mathbf{V}}_{\boldsymbol{c}}^{-1} \boldsymbol{1} - (\boldsymbol{1}^{\top} \widehat{\mathbf{V}}_{\boldsymbol{c}}^{-1} \widehat{\boldsymbol{B}}_{\boldsymbol{c}}) \widehat{\mathbf{V}}_{\boldsymbol{c}}^{-1}\widehat{\boldsymbol{B}}_{\boldsymbol{c}}}{(1+\widehat{\boldsymbol{B}}_{\boldsymbol{c}}^{\top} \widehat{\mathbf{V}}_{\boldsymbol{c}}^{-1}\widehat{\boldsymbol{B}}_{\boldsymbol{c}}) (\boldsymbol{1}^{\top} \widehat{\mathbf{V}}_{\boldsymbol{c}}^{-1} \boldsymbol{1}) - (\boldsymbol{1}^{\top} \widehat{\mathbf{V}}_{\boldsymbol{c}}^{-1} \widehat{\boldsymbol{B}}_{\boldsymbol{c}})^2}.
\]
Next, we provide the asymptotic properties of the composite pooled estimators based on these estimated optimal values of weights. 
\begin{Coro}
\label{coro:pool_mfixed_kinfty_general}
Work under the conditions of Theorem~\ref{theo:pool_mfixed_kinfty_general} with $\boldsymbol{\gamma}=\gamma \boldsymbol{1}$, $\rho_j<0$, and $A_j(t)=\gamma \beta_j t^{\rho_j}$ for all $j$. Assume that the matrix $\mathbf{V}$ is positive definite (hence $\mathbf{V}_{\boldsymbol{c}}$ is). Assume further that, for all $j\in\{1,\ldots,m\}$, $\widehat{\beta}_j$ is a consistent estimator of $\beta_j$ and $(\widehat{\rho}_j-\rho_j) \log n_j=\operatorname{o}_{\mathbb{P}}(1)$. Then
\begin{align*}
&\sqrt{k} (\widehat{\gamma}_n(\widehat{\boldsymbol{\omega}}_n^{(\mathrm{Var})})-\gamma) \stackrel{d}{\longrightarrow} \mathcal{N}\left( \dfrac{\boldsymbol{1}^{\top} \mathbf{V}_{\boldsymbol{c}}^{-1} \boldsymbol{B}_{\boldsymbol{c}}}{\boldsymbol{1}^{\top} \mathbf{V}_{\boldsymbol{c}}^{-1} \boldsymbol{1}}, \dfrac{1}{\boldsymbol{1}^{\top} \mathbf{V}_{\boldsymbol{c}}^{-1} \boldsymbol{1}} \right), \mbox{ and } \\
&\sqrt{k} (\widehat{\gamma}_n(\widehat{\boldsymbol{\omega}}_n^{(\mathrm{AMSE})})-\gamma) \stackrel{d}{\longrightarrow} \mathcal{N}\left( (\boldsymbol{\omega}^{(\mathrm{AMSE})})^{\top} \boldsymbol{B}_{\boldsymbol{c}}, (\boldsymbol{\omega}^{(\mathrm{AMSE})})^{\top} \mathbf{V}_{\boldsymbol{c}} \boldsymbol{\omega}^{(\mathrm{AMSE})} \right), \mbox{ where } \\
&((\boldsymbol{\omega}^{(\mathrm{AMSE})})^{\top} \boldsymbol{B}_{\boldsymbol{c}})^2 + (\boldsymbol{\omega}^{(\mathrm{AMSE})})^{\top} \mathbf{V}_{\boldsymbol{c}} \boldsymbol{\omega}^{(\mathrm{AMSE})} = \frac{1+\boldsymbol{B}_{\boldsymbol{c}}^{\top} \mathbf{V}_{\boldsymbol{c}}^{-1} \boldsymbol{B}_{\boldsymbol{c}}}{(1+\boldsymbol{B}_{\boldsymbol{c}}^{\top} \mathbf{V}_{\boldsymbol{c}}^{-1} \boldsymbol{B}_{\boldsymbol{c}}) (\boldsymbol{1}^{\top} \mathbf{V}_{\boldsymbol{c}}^{-1} \boldsymbol{1}) - (\boldsymbol{1}^{\top} \mathbf{V}_{\boldsymbol{c}}^{-1} \boldsymbol{B}_{\boldsymbol{c}})^2}.
\end{align*}
\end{Coro}
\begin{Rmk}[On the variance- and AMSE-optimal choices]
\label{rmk:varopt_convex}
This is the first work to implement the idea of AMSE-optimal weights $\boldsymbol{\omega}^{(\mathrm{AMSE})}$ which should be favored in practice when the estimation bias is large enough, especially in the case of unequal sample fractions $k_j/n_j$, as demonstrated below in Section~\ref{sec:distrinf:var_AMSE} in the distributed inference framework. The solution using variance-optimal weights $\boldsymbol{\omega}^{(\mathrm{Var})}$ was actually suggested by~\cite{kinfrilil2016}, but their analog of Theorem~\ref{theo:pool_mfixed_kinfty_general} requires much stronger technical conditions. Besides, the use of their estimator for $\boldsymbol{\omega}^{(\mathrm{Var})}$ lacks a theoretical justification similar to Corollary~\ref{coro:pool_mfixed_kinfty_general}. A different, variance-optimal convex combination (with nonnegative weights) is advocated in~\cite{demcle2016}, also under much stronger technical conditions. Their solution does not coincide with ours in general, although this is not so clear-cut in the specific case $c_j=1$ of identical effective sample sizes. In this situation, 
at least when $m\leq 3$, the variance-optimal set of weights we propose is in fact a convex combination. For $m=2$, a direct calculation provides $\boldsymbol{\omega}^{(\mathrm{Var})}=(1/2,1/2)^{\top}$ irrespective of 
$\mathbf{V}_{\boldsymbol{1}} = \mathbf{V}$, corresponding to the naive average. For $m=3$, the discussion is more complex and involves the identification of tail correlation matrices with a convex polytope that is a proper subset of the {\it elliptope}, a Riemannian quotient manifold representing the set of standard correlation matrices. As a consequence, pooling together $m\leq 3$ Hill estimators with equal effective sample size can never 
outperform, in terms of asymptotic variance, a Hill estimator built from a pooled sample of independent data of equivalent total size. This would not necessarily be the case in general pooling problems, 
{\it e.g.}, for positively correlated sample means. What can happen in the case $m>3$ remains unclear. 
\end{Rmk}
We conclude this section by discussing bias-reduced versions of the variance-optimal and AMSE-optimal pooled estimators, defined as
\begin{align*}
\overline{\gamma}_n(\widehat{\boldsymbol{\omega}}_n^{(\mathrm{Var})}) &= \widehat{\gamma}_n(\widehat{\boldsymbol{\omega}}_n^{(\mathrm{Var})}) - \frac{1}{\sqrt{k}}(\widehat{\boldsymbol{\omega}}_n^{(\mathrm{Var})})^{\top} \widehat{\boldsymbol{B}}_{\boldsymbol{c}} \\
\mbox{and } \overline{\gamma}_n(\widehat{\boldsymbol{\omega}}_n^{(\mathrm{AMSE})}) &= \widehat{\gamma}_n(\widehat{\boldsymbol{\omega}}_n^{(\mathrm{AMSE})}) - \frac{1}{\sqrt{k}}(\widehat{\boldsymbol{\omega}}_n^{(\mathrm{AMSE})})^{\top} \widehat{\boldsymbol{B}}_{\boldsymbol{c}}.
\end{align*}
\begin{Coro}
\label{coro:pool_mfixed_kinfty_general_BR}
Under the conditions of Corollary~\ref{coro:pool_mfixed_kinfty_general},
$\sqrt{k} (\overline{\gamma}_n(\widehat{\boldsymbol{\omega}}_n^{(\mathrm{Var})})-\gamma) \stackrel{d}{\longrightarrow} \mathcal{N}( 0, 1/(\boldsymbol{1}^{\top} \mathbf{V}_{\boldsymbol{c}}^{-1} \boldsymbol{1}) )$ and
\begin{align*}
    & \sqrt{k} (\overline{\gamma}_n(\widehat{\boldsymbol{\omega}}_n^{(\mathrm{AMSE})})-\gamma) \\ &\stackrel{d}{\longrightarrow} \mathcal{N}\left( 0, \frac{(1+\boldsymbol{B}_{\boldsymbol{c}}^{\top} \mathbf{V}_{\boldsymbol{c}}^{-1}\boldsymbol{B}_{\boldsymbol{c}})^2 (\boldsymbol{1}^{\top} \mathbf{V}_{\boldsymbol{c}}^{-1} \boldsymbol{1}) - (2+\boldsymbol{B}_{\boldsymbol{c}}^{\top} \mathbf{V}_{\boldsymbol{c}}^{-1}\boldsymbol{B}_{\boldsymbol{c}}) (\boldsymbol{1}^{\top} \mathbf{V}_{\boldsymbol{c}}^{-1} \boldsymbol{B}_{\boldsymbol{c}})^2}{[(1+\boldsymbol{B}_{\boldsymbol{c}}^{\top} \mathbf{V}_{\boldsymbol{c}}^{-1}\boldsymbol{B}_{\boldsymbol{c}}) (\boldsymbol{1}^{\top} \mathbf{V}_{\boldsymbol{c}}^{-1} \boldsymbol{1}) - (\boldsymbol{1}^{\top} \mathbf{V}_{\boldsymbol{c}}^{-1} \boldsymbol{B}_{\boldsymbol{c}})^2]^2} \right).
\end{align*}
\end{Coro}
If the marginal distributions are equal across samples, then one can improve the estimation of the weights by pooling the second-order parameter estimators themselves. We shall explore this possibility in Section~\ref{sec:distrinf:var_AMSE}.
\subsection{Weighted geometric pooling of extreme quantile estimators}
\label{sec:main:extquant}

We turn to extreme quantile estimation for a very small exceedance probability $p=p(n)\to 0$ as $n\to\infty$. In each sample, we estimate $q_j(1-p)$ by the extrapolated Weissman estimator of~\cite{wei1978}:
\[
\widehat{q}_j^{\star}(1-p|k_j) = \left( \frac{k_j}{n_j p} \right)^{\widehat{\gamma}_j(k_j)} X_{n_j-k_j:n_j,j},
\]
where $k_j$ satisfies $k_j/(n_j p)\to\infty$. The typical case is when $np$ is bounded, reflecting the challenging problem of estimating quantiles in the far tail where only a few or no observations are available in the neighborhood of $q_j(1-p)$. When the samples are believed or known to have the same tail index $\gamma$, it is natural to harness the power of pooling by substituting the weighted estimator $\widehat{\gamma}_n(\boldsymbol{\omega})$ in place of the individual estimator $\widehat{\gamma}_j(k_j)$, to get 
\[
\widehat{q}_j^{\star}(1-p|k_j,\boldsymbol{\omega}) = \left( \frac{k_j}{n_j p} \right)^{\widehat{\gamma}_n(\boldsymbol{\omega})} X_{n_j-k_j:n_j,j}. 
\]
This improves on the traditional Weissman estimator in each sample by borrowing statistical strength across all samples. Going one step further, the marginal quantile estimators themselves can be pooled to gain more insight when the datasets have the same extreme quantiles, or equivalently, if one assumes that
\vskip1ex
\noindent
$(\mathcal{H})$ For any $j,\ell\in \{1,\ldots,m\}$ with $j\neq \ell$, we have $U_j(t)/U_{\ell}(t)\to 1$ as $t\to \infty$. 
\vskip1ex
\noindent
The validity of this assumption can be tested by applying the test of~\cite{padstu2021} described in their Sections~3.3 and~5.
Under assumption $(\mathcal{H})$, the quantiles $q_j(1-p)$ are all asymptotically equivalent and hence the estimators $\widehat{q}_j^{\star}(1-p|k_j)$ all estimate the same quantity. A straightforward way of combining these individual extreme quantile estimators would be to take again a weighted sum, as $\overline{q}_n^{\star}(1-p|\boldsymbol{\omega})=\sum_{j=1}^m \omega_j \widehat{q}_j^{\star}(1-p|k_j)$. This is most efficient when estimating central quantiles, see for example~\cite{knibas2003}. However, it is no longer the best possible solution when it comes to pooling the Weissman estimators, because the use of geometric weighted sums better suits the multiplicative and power structure of these extrapolated estimators (see Sections~\ref{sec:fin:simul:genpool} for numerical evidence). The crucial point to note here is that the log-Weissman quantile estimator can be rewritten as 
\[
\log \widehat{q}_j^{\star}(1-p|k_j) = \log \left( \frac{k}{n p} \right) \widehat{\gamma}_j(k_j) + \left[ \log \left( \frac{k_j}{k} \right) - \log \left( \frac{n_j}{n} \right) \right] \widehat{\gamma}_j(k_j) +  \log X_{n_j-k_j:n_j,j}. 
\]
In the first term, which dominates since $k/(np)\to\infty$, the $\gamma$ estimator now appears on the standard scale. This suggests the use of
$
\log \widehat{q}_n^{\star}(1-p|\boldsymbol{\omega}) = \sum_{j=1}^m \omega_j \log \widehat{q}_j^{\star}(1-p|k_j)
$
or, in other words, 
\[
\widehat{q}_n^{\star}(1-p|\boldsymbol{\omega}) = \prod_{j=1}^m [\widehat{q}_j^{\star}(1-p|k_j)]^{\omega_j} = \prod_{j=1}^m \left[ \left( \frac{k_j}{n_j p} \right)^{\widehat{\gamma}_j(k_j)} X_{n_j-k_j:n_j,j} \right]^{\omega_j}.
\]
This estimator is a weighted geometric (rather than arithmetic) mean of the $\widehat{q}_j^{\star}(1-p|k_j)$. We conclude this discussion by deriving the asymptotic normality of $\widehat{q}_j^{\star}(1-p|k_j,\boldsymbol{\omega})$ and $\widehat{q}_n^{\star}(1-p|\boldsymbol{\omega})$. 
\begin{Theo}
\label{theo:pool_mfixed_kinfty_general_W}
Work under the conditions and with the notation of Theorem~\ref{theo:pool_mfixed_kinfty_general} with $\boldsymbol{\gamma}=\gamma \boldsymbol{1}$ and $\rho_j<0$ for all $j\in \{1,\ldots,m\}$. Pick $p=p(n)\to 0$ such that $k/(n p)\to \infty$ and $\sqrt{k}/\log(k/(n p))\to \infty$ as $n\to\infty$. Let $\boldsymbol{\omega}, \widehat{\boldsymbol{\omega}}_n$ be such that $\widehat{\boldsymbol{\omega}}_n^{\top} \boldsymbol{1}=1$ and $\widehat{\boldsymbol{\omega}}_n\stackrel{\mathbb{P}}{\longrightarrow} \boldsymbol{\omega}$. Then, for any $j$, 
\[
\frac{\sqrt{k}}{\log(k_j/(n_j p))} \left( \frac{\widehat{q}_j^{\star}(1-p|k_j,\widehat{\boldsymbol{\omega}}_n)}{q_j(1-p)} - 1 \right) = \sqrt{k} (\widehat{\gamma}_n(\boldsymbol{\omega}) - \gamma) + \operatorname{o}_{\mathbb{P}}(1) \stackrel{d}{\longrightarrow} \mathcal{N}\left( \boldsymbol{\omega}^{\top} \boldsymbol{B}_{\boldsymbol{c}}, \boldsymbol{\omega}^{\top} \mathbf{V}_{\boldsymbol{c}} \boldsymbol{\omega} \right).
\]
If moreover assumption $(\mathcal{H})$ holds then, for any $j$, 
\[
    \frac{\sqrt{k}}{\log(k/(n p))} \left( \frac{\widehat{q}_n^{\star}(1-p|\widehat{\boldsymbol{\omega}}_n)}{q_j(1-p)} - 1 \right) = \sqrt{k} (\widehat{\gamma}_n(\boldsymbol{\omega}) - \gamma) + \operatorname{o}_{\mathbb{P}}(1)  \stackrel{d}{\longrightarrow} \mathcal{N}\left( \boldsymbol{\omega}^{\top} \boldsymbol{B}_{\boldsymbol{c}}, \boldsymbol{\omega}^{\top} \mathbf{V}_{\boldsymbol{c}} \boldsymbol{\omega} \right).
\]
\end{Theo}
An analogue of Theorem~\ref{theo:pool_mfixed_kinfty_general_W} is feasible for optimally-pooled extreme quantile estimation where $\widehat{\boldsymbol{\omega}}_n\in \{ \widehat{\boldsymbol{\omega}}_n^{(\mathrm{Var})}, \widehat{\boldsymbol{\omega}}_n^{(\mathrm{AMSE})} \}$, since the asymptotic distribution of $\widehat{q}_n^{\star}(1-p|\widehat{\boldsymbol{\omega}}_n)$ is governed by that of $\widehat{\gamma}_n(\boldsymbol{\omega})$. Similar results can also be established when $\widehat{\gamma}_n(\boldsymbol{\omega})$ is replaced by the bias-reduced versions $\overline{\gamma}_n(\widehat{\boldsymbol{\omega}}_n^{(\mathrm{Var})})$ and $\overline{\gamma}_n(\widehat{\boldsymbol{\omega}}_n^{(\mathrm{AMSE})})$, so they are omitted. 
\subsection{Inference using pooled extreme value estimators}
\label{sec:main:inference}

Unless there are strong reasons to believe in the equality of tail indices (as is the case in the distributed inference setup of Section~\ref{sec:distrinf}), it is crucial to justify this assumption by performing a statistical test before applying our pooled estimators. To do so, we briefly present here an approach motivated by testing for nested models. Suppose that $\boldsymbol{Z}$ is an $m-$dimensional Gaussian random vector with mean $\boldsymbol{\mu}$ and known positive definite covariance matrix $\boldsymbol{V}$, and consider the testing problem of
$
M_0: \mu_{1}=\cdots=\mu_{m}=\mu
$
versus 
$
M_1: \exists(j,\ell) \mbox{ with } j\neq \ell \mbox{ such that } \mu_j\neq \mu_{\ell},
$
based on $\boldsymbol{Z}$. The log-likelihood ratio deviance statistic for testing the validity of model $M_0$ is $\Lambda = (\boldsymbol{Z} - \widehat{\mu} \boldsymbol{1})^{\top} \mathbf{V}^{-1} (\boldsymbol{Z} - \widehat{\mu} \boldsymbol{1})$, with $\widehat{\mu} = (\boldsymbol{1}^{\top} \mathbf{V}^{-1} \boldsymbol{Z})/(\boldsymbol{1}^{\top} \mathbf{V}^{-1} \boldsymbol{1})$. In model $M_0$, the statistic $\Lambda$ has a chi-square distribution with $m-1$ degrees of freedom. In our context, under the assumptions of Theorem~\ref{theo:pool_mfixed_kinfty_general} and if all the $\lambda_j$ are 0 (see Remark~\ref{rmk:inference} below for more discussion on this assumption), one has 
\begin{align*}
    \sqrt{k}(\widehat{\boldsymbol{\gamma}}_n - \boldsymbol{\gamma}) &= (\sqrt{k} (\widehat{\gamma}_1(k_1)-\gamma_1),\ldots,\sqrt{k}(\widehat{\gamma}_m(k_m)-\gamma_m))^{\top} \stackrel{d}{\longrightarrow} \mathcal{N}( \boldsymbol{0}, \mathbf{V}_{\boldsymbol{c}}) \\
    \mbox{with } & [\mathbf{V}_{\boldsymbol{c}}]_{j,\ell} = \left( \sum_{i=1}^m c_i^{-1} \right) \begin{cases} \gamma_j^2 c_j  & \mbox{if } j=\ell, \\[5pt] \gamma_j \gamma_{\ell} c_j R_{j,\ell}(c_{\ell}/c_j,b_{\ell}/b_j) & \mbox{if } j<\ell. \end{cases} 
\end{align*}
This 
can be formulated by the approximation~$\widehat{\boldsymbol{\gamma}}_n\stackrel{d}{\approx} \mathcal{N}( \boldsymbol{\gamma}, k^{-1} \mathbf{V}_{\boldsymbol{c}} )$ of the distribution of $\widehat{\boldsymbol{\gamma}}_n$ by the $\mathcal{N}( \boldsymbol{\gamma}, k^{-1} \mathbf{V}_{\boldsymbol{c}} )$ distribution. Given the estimator 
\[
[\overline{\mathbf{V}}_{\boldsymbol{c}}]_{j,\ell} = 
k \begin{cases} \widehat{\gamma}_j^2(k_j) \dfrac{1}{k_j}  & \mbox{if } j=\ell, \\[10pt] \widehat{\gamma}_j(k_j) \widehat{\gamma}_{\ell}(k_{\ell}) \dfrac{1}{k_j} \widehat{R}_{j,\ell}(k_j/k_{\ell},n_j/n_{\ell}) & \mbox{if } j\neq\ell, \end{cases} 
\]
one can obtain a deviance statistic for testing $H_0: \boldsymbol{\gamma}=\gamma\boldsymbol{1}$ versus $H_1:\boldsymbol{\gamma}\neq \gamma\boldsymbol{1}$ as 
\[
\Lambda_n = k (\widehat{\boldsymbol{\gamma}}_n - \overline{\mu}_n \boldsymbol{1})^{\top} \overline{\mathbf{V}}_{\boldsymbol{c}}^{-1} (\widehat{\boldsymbol{\gamma}}_n - \overline{\mu}_n \boldsymbol{1}), \mbox{ with } \overline{\mu}_n = \frac{\boldsymbol{1}^{\top} \overline{\mathbf{V}}_{\boldsymbol{c}}^{-1} \widehat{\boldsymbol{\gamma}}_n}{\boldsymbol{1}^{\top} \overline{\mathbf{V}}_{\boldsymbol{c}}^{-1} \boldsymbol{1}} = (\overline{\boldsymbol{\omega}}_n^{(\mathrm{Var})})^{\top} \widehat{\boldsymbol{\gamma}}_n = \widehat{\gamma}_n(\overline{\boldsymbol{\omega}}_n^{(\mathrm{Var})}). 
\]
Therefore the test statistic $\Lambda_n$ compares the vector of estimates $\widehat{\boldsymbol{\gamma}}_n$ with an estimate of the variance-optimal pooled estimator on a scale adapted to the amount of dependence existing between the extremes of the vector $(X_1,\ldots,X_m)^{\top}$. A somewhat different proposal, not motivated by a likelihood ratio test in nested models, is outlined in~\cite{kinfrilil2016}. 

Our testing procedure, of asymptotic size $\alpha$, is to reject $H_0$ if $\Lambda_n>\chi_{m-1,1-\alpha}^2$, where $\chi_{m-1,1-\alpha}^2$ is the $(1-\alpha)$th quantile of the chi-square distribution with $m-1$ degrees of freedom. The next corollary establishes the consistency of this test and gives a symmetric asymptotic confidence interval for the common tail index $\gamma$ under $H_0$.
\begin{Coro}
\label{coro:pool_testing}
Under the conditions of Corollary~\ref{coro:pool_mfixed_kinfty_general} and the assumption that $\lambda_j=0$ for all $j$, we have $\mathbb{P}(\Lambda_n>\chi^2_{m-1,1-\alpha}) \to \alpha$ under $H_0$, and $\Lambda_n\stackrel{\mathbb{P}}{\longrightarrow} +\infty$ under $H_1$, as $n\to\infty$. Moreover, under $H_0$, if $\widehat{\boldsymbol{\omega}}_n^{\top} \boldsymbol{1}=1$ with $\widehat{\boldsymbol{\omega}}_n\stackrel{\mathbb{P}}{\longrightarrow} \boldsymbol{\omega}$ then, for any $\alpha\in (0,1)$, 
\[
\lim_{n\to\infty} \mathbb{P}\left( \gamma \in \left[ \widehat{\gamma}_n(\widehat{\boldsymbol{\omega}}_n) \pm z_{1-\alpha/2} \sqrt{(\widehat{\boldsymbol{\omega}}_n^{\top} \widehat{\mathbf{V}}_{\boldsymbol{c}} \widehat{\boldsymbol{\omega}}_n)/k} \right] \right) = 1-\alpha
\]
where $z_{1-\alpha/2}$ is the $(1-\alpha/2)$th quantile of the standard Gaussian distribution. [In this asymptotic confidence interval, $\widehat{\mathbf{V}}_{\boldsymbol{c}}$ is calculated as described in Section~\ref{sec:main:optimal_EVI}.] 
\end{Coro}
\begin{Rmk}[With asymptotic independence across subsamples]
\label{rmk:asyindep}
An important subcase in practice is when pairs of data points taken from two different subsamples are asymptotically independent. This covers, for example, the distributed inference situation which will be discussed in Section~\ref{sec:distrinf}. In this case, all tail copulae $R_{j,\ell}$ are identically zero, so one can estimate $\mathbf{V}_{\boldsymbol{c}}$ with $k\operatorname{diag}( \widehat{\gamma}_1^2(k_1)/k_1,\ldots,\widehat{\gamma}_m^2(k_m)/k_m)$. The test statistic $\Lambda_n$ becomes 
\[
\Lambda_n=\sum_{j=1}^m k_j \frac{(\widehat{\gamma}_j(k_j) - \widehat{\gamma}_n(\overline{\boldsymbol{\omega}}_n^{(\mathrm{Var})}))^2}{\widehat{\gamma}_j^2(k_j)}.
\]
This has the familiar look of a Pearson goodness-of-fit statistic, with the weight $k_j$ adjusting for the different rates of convergence of the $\widehat{\gamma}_j(k_j)$. If all the $k_j$ are equal, then 
\[
\Lambda_n=\frac{k}{m} \sum_{j=1}^m \left( \frac{\widehat{\gamma}_n(\overline{\boldsymbol{\omega}}_n^{(\mathrm{Var})})}{\widehat{\gamma}_j(k_j)} - 1 \right)^2. 
\]
In this setup, our proposed statistic $\Lambda_n$ bears some similarity with a test statistic studied in~\cite{einfanzho2020} in the context of testing for the validity of a multivariate regular variation model that assumes equality of tail indices across marginal distributions.
\end{Rmk}
\begin{Rmk}[Inference and bias correction]
\label{rmk:inference}
Typically, assuming $\lambda_j=0$ to omit the asymptotic bias terms is sensible as long as the second-order parameters $\rho_j$ remain reasonably far away from 0. Based on finite-sample experiments with a total sample size $n=1{,}000$, marginal Burr distributions and $2\leq m\leq 5$ with both balanced and unbalanced samples, the confidence interval provided seems to perform very well when $|\rho_j|>3/4$. 
Estimating the bias terms then is in fact detrimental, because of increased variability of the resulting interval estimator that is not accounted for in the estimated variance.
\end{Rmk}
\begin{Rmk}[Tail homogeneity and tail homoskedasticity]
\label{rmk:hetero}
When all the parameters $\rho_j$ are negative, as in Corollary~\ref{coro:pool_testing}, one has $t^{-\gamma_j} U_j(t)\to C_j\in (0,\infty)$ as $t\to\infty$, see the equation below Equation~(2.3.23) in~\cite{haafer2006}. 
Testing $H_0: \boldsymbol{\gamma}=\gamma\boldsymbol{1}$ versus $H_1:\boldsymbol{\gamma}\neq \gamma\boldsymbol{1}$ is then exactly equivalent to testing
\[
\begin{array}{rl}
 & H'_0 : \forall j,\ell \in \{1,\ldots,d\}, \ \lim_{\tau\uparrow 1} q_j(\tau)/q_{\ell}(\tau)\in (0,\infty), \\[3pt]
\mbox{versus } & H'_1 : \exists j,\ell \in \{1,\ldots,d\} \mbox{ with } j\neq \ell  \mbox{ and } \lim_{\tau\uparrow 1} q_j(\tau)/q_{\ell}(\tau) \in \{ 0,\infty \}.
\end{array} 
\]
The testing procedure based on the statistic $\Lambda_n$ is therefore, under very mild conditions, exactly a test for asymptotic proportionality of marginal extreme quantiles. It can thus be used to detect tail homogeneity (equal tail indices and therefore asymptotically proportional tail quantiles) as opposed to tail heterogeneity (one marginal distribution having a heavier tail than the others). We discuss below the testing of the stronger property when all limits are equal to 1 in $H'_0$, corresponding to the asymptotic equivalence of extreme quantiles $(\mathcal{H})$, and referred to as tail homoskedasticity.
\end{Rmk}
Testing for tail homoskedasticity can be done directly using the Weissman estimators $\widehat{q}_j^{\star}(1-p|k_j)$. Set $\boldsymbol{Z}_n(p)=\log \widehat{\boldsymbol{q}}_n^{\star}(1-p)=(\log \widehat{q}_1^{\star}(1-p|k_1),\ldots,\log \widehat{q}_m^{\star}(1-p|k_m))$ and 
\[
L_n(p) = \frac{k}{\log^2(k/(np))} \left( \boldsymbol{Z}_n(p) - \frac{\boldsymbol{1}^{\top} \overline{\mathbf{V}}_{\boldsymbol{c}}^{-1} \boldsymbol{Z}_n(p)}{\boldsymbol{1}^{\top} \overline{\mathbf{V}}_{\boldsymbol{c}}^{-1} \boldsymbol{1}} \boldsymbol{1} \right)^{\top} \overline{\mathbf{V}}_{\boldsymbol{c}}^{-1} \left( \boldsymbol{Z}_n(p) - \frac{\boldsymbol{1}^{\top} \overline{\mathbf{V}}_{\boldsymbol{c}}^{-1} \boldsymbol{Z}_n(p)}{\boldsymbol{1}^{\top} \overline{\mathbf{V}}_{\boldsymbol{c}}^{-1} \boldsymbol{1}} \boldsymbol{1} \right). 
\]
A testing procedure of asymptotic size $\alpha$
of $(\mathcal{H})$ versus
\[
(\mathcal{H}'): \exists j,\ell \in \{1,\ldots,d\} \mbox{ with } j\neq \ell \mbox{ and } \lim_{\tau\uparrow 1} q_j(\tau)/q_{\ell}(\tau) \neq 1,
\]
is to reject $(\mathcal{H})$ if $L_n(p)>\chi_{m-1,1-\alpha}^2$. We establish this rigorously in our next result.
\begin{Coro}
\label{coro:pool_testing_quantile}
Under the conditions of Corollary~\ref{coro:pool_testing} and $\rho_j<0$ for all $j$, if $p=p(n)\to 0$ is such that $k/(n p)\to \infty$ and $\sqrt{k}/\log(k/(n p))\to \infty$ as $n\to\infty$, then we have $\mathbb{P}(L_n(p)>\chi^2_{m-1,1-\alpha}) \to \alpha$ under $(\mathcal{H})$, and $L_n(p)\stackrel{\mathbb{P}}{\longrightarrow} +\infty$ as $n\to\infty$ under $(\mathcal{H}')$. Moreover, under $(\mathcal{H})$, if $\widehat{\boldsymbol{\omega}}_n^{\top} \boldsymbol{1}=1$ with $\widehat{\boldsymbol{\omega}}_n\stackrel{\mathbb{P}}{\longrightarrow} \boldsymbol{\omega}$ then, for any $\alpha\in (0,1)$ and $j\in \{ 1,\ldots,m\}$, 
\[
\lim_{n\to\infty} \mathbb{P}\left( q_j(1-p) \in \left[ \widehat{q}_n^{\star}(1-p|\widehat{\boldsymbol{\omega}}_n) \exp\left( \pm z_{1-\alpha/2} \log\left[\frac{k}{np}\right] \sqrt{\frac{\widehat{\boldsymbol{\omega}}_n^{\top} \widehat{\mathbf{V}}_{\boldsymbol{c}} \widehat{\boldsymbol{\omega}}_n}{k}} \right) \right] \right) = 1-\alpha.
\]
\end{Coro}
The present test 
is more general than the one suggested in Section~3.3 of~\cite{padstu2021}, for which fairly strong integrability assumptions on the $X_j$ are unavoidable.
The use of the log-scale is equivalent in theory to the relative scale employed in Theorem~\ref{theo:pool_mfixed_kinfty_general_W}, but it tends to provide more accurate asymptotic confidence intervals for extreme quantiles, as indicated by {\it e.g.}~\cite{dre2003}.

\section{The framework of distributed inference}
\label{sec:distrinf}

Our general pooling theory naturally applies to the context of distributed inference. A restriction in this 
framework 
is that the~$n$ individual data points cannot be processed with a standalone machine and very restricted or no communication is allowed between the~$m$ machines. In particular, the end user operating from a central machine only has access to limited information, such as the subsample estimates and associated~$n_j$ and~$k_j$, which is not sufficient for estimating the tail dependence structure between the different subsamples. We thus make the assumption that the data within and across machines are independent, that is, the $X_{i,j}$ are i.i.d.~for $1\leq i\leq n_j$ and $1\leq j\leq m$, with a common distribution satisfying the second-order condition $\mathcal{C}_2(\gamma,\rho,A)$.
\subsection{Distributed estimation of the tail index}
\label{sec:distrinf:pooling}

When the data are i.i.d., an obvious benchmark for the 
distributed estimator $\widehat{\gamma}_n(\boldsymbol{\omega})$ is the Hill estimator 
based on the 
unfeasible combination (due to computational or storage difficulties) of subsamples $\{ X_i, 1\leq i\leq n\} = \{ X_{i,j}, 1\leq j\leq m, 1\leq i\leq n_j \}$ with effective sample size $k=\sum_{j=1}^m k_j$, that is,
$
\widehat{\gamma}_n^{(\mathrm{Hill})}(k) = k^{-1} \sum_{i=1}^k \log (X_{n-i+1:n}/X_{n-k:n})
$
where $X_{1:n}\leq X_{2:n}\leq \cdots\leq X_{n:n}$ are the order statistics of the random variables $X_i$. Assume, as in Section~\ref{sec:main:EVI}, that the $n_j$ and the $k_j$ are asymptotically proportional but possibly unbalanced, {\it i.e.}~$n_1/n_j\to b_j\in (0,\infty)$ and $k_1/k_j\to c_j\in (0,\infty)$. Then, since $A$ is regularly varying with index $\rho$, it is straightforward to show that 
the existence of 
$\lambda_j=\lim_{n\to\infty} \sqrt{k_j} A(n_j/k_j)$, for any $j$, is equivalent to the existence of $\lambda=\lim_{n\to\infty} \sqrt{k} A(n/k)$, and that
\[
\lambda_j = c_j^{\rho-1/2} b_j^{-\rho} \lambda_1 = c_j^{\rho-1/2} b_j^{-\rho} \left( \sum_{j=1}^m \frac{1}{c_j} \right)^{\rho-1/2} \left( \sum_{j=1}^m \frac{1}{b_j} \right)^{-\rho} \lambda. 
\]
Hence, we have the following corollary of Theorem~\ref{theo:pool_mfixed_kinfty_general}. 
\begin{Coro}
\label{coro:pool_mfixed_kinfty_balanced}
Assume that condition $\mathcal{C}_2(\gamma,\rho,A)$ holds. Suppose that $n_1/n_j\to b_j\in (0,\infty)$ and $k_1/k_j\to c_j\in (0,\infty)$ (with then $b_1,c_1=1$) for any $j\in \{1,\ldots,m\}$, and then that $k\to\infty$ with $k/n\to 0$ and $\sqrt{k} A(n/k)\to \lambda\in\mathbb{R}$, as $n\to\infty$. Let $\boldsymbol{\omega}=(\omega_1,\ldots,\omega_m)^{\top}$ be such that $\boldsymbol{\omega}^{\top} \boldsymbol{1}=1$. Then 
\[
\sqrt{k} (\widehat{\gamma}_n(\boldsymbol{\omega})-\gamma) \stackrel{d}{\longrightarrow} \mathcal{N}\left( \frac{\lambda}{1-\rho} \sum_{j=1}^m d_j^{\rho} \omega_j, \gamma^2 \sum_{j=1}^m \frac{1}{c_j} \sum_{j=1}^m c_j \omega_j^2 \right)
\]
where $d_j=(c_j/b_j) \times (\sum_{i=1}^m c_i^{-1})/(\sum_{i=1}^m b_i^{-1})$. If $\widehat{\boldsymbol{\omega}}_n^{\top} \boldsymbol{1}=1$ and $\widehat{\boldsymbol{\omega}}_n\stackrel{\mathbb{P}}{\longrightarrow} \boldsymbol{\omega}$ then $\sqrt{k}(\widehat{\gamma}_n(\widehat{\boldsymbol{\omega}}_n)-\widehat{\gamma}_n(\boldsymbol{\omega}))=\operatorname{o}_{\mathbb{P}}(1)$ and so the above convergence 
remains valid for $\widehat{\gamma}_n(\widehat{\boldsymbol{\omega}}_n)$.
\end{Coro}
With this 
result, we discuss next in full detail the performance 
of the variance- and AMSE-optimal distributed estimators relative to the benchmark Hill estimator $\widehat{\gamma}_n^{(\mathrm{Hill})}(k)$, which satisfies the weak convergence
$
\sqrt{k} (\widehat{\gamma}_n^{(\mathrm{Hill})}(k)-\gamma) \stackrel{d}{\longrightarrow} \mathcal{N}( \lambda/(1-\rho), \gamma^2)
$ under the assumptions of Corollary~\ref{coro:pool_mfixed_kinfty_balanced} (see Theorem~3.2.5, p.74 
in~\cite{haafer2006}).

\subsection{Variance-optimal and AMSE-optimal combinations}
\label{sec:distrinf:var_AMSE}

It is immediate that the variance-optimal weights are $\omega_j^{(\mathrm{Var})}=(\sum_{i=1}^m c_i^{-1} )^{-1} c_j^{-1}$ for all $j$. Estimating $1/c_i$ by $k_i/k_1$ leads to the estimated weight vector $\widetilde{\boldsymbol{\omega}}_n^{(\mathrm{Var})} = (k_1/k,\ldots,k_m/k)$. The variance-optimal distributed estimator is therefore a very convenient convex combination that has a much simpler expression than in the general setup of Section~\ref{sec:main}, due to the constraint of data independence across machines. 
Calculating this estimator only requires reporting each individual $k_j$ and $\widehat\gamma_j(k_j)$ to the central machine. Being a nonrandom convex combination, this estimator is immune to the instability issues caused by possible inaccurate estimation of optimal weights.nThe following result gives its asymptotic distribution.
\begin{Coro}
\label{coro:poolopt_mfixed_kinfty_bal}
Under the conditions and notation of Corollary~\ref{coro:pool_mfixed_kinfty_balanced}, 
\[ 
\sqrt{k} (\widehat{\gamma}_n(\widetilde{\boldsymbol{\omega}}_n^{(\mathrm{Var})})-\gamma) \stackrel{d}{\longrightarrow} \mathcal{N}\left( \frac{\lambda}{1-\rho} \left( \sum_{j=1}^m \frac{1}{c_j} \right)^{-1} \left( \sum_{j=1}^m \frac{d_j^{\rho}}{c_j} \right), \gamma^2 \right).
\]
\end{Coro}
\begin{Rmk}[Asymptotic bias comparison] 
\label{rmk:bias} 
The unfeasible Hill estimator
has asymptotic bias $\mu^{(\mathrm{Hill})}=\lambda/(1-\rho)$. 
If $\lambda\neq 0$ 
and $\mu^{(\mathrm{Var})}$ denotes the asymptotic bias of $\widehat{\gamma}_n(\widetilde{\boldsymbol{\omega}}_n^{(\mathrm{Var})})$, then
clearly $\mu^{(\mathrm{Hill})}=\mu^{(\mathrm{Var})}$ when $\rho=0$. Otherwise, the H\"older inequality for the conjugate exponents $p=-(1-\rho)/\rho$ and $q=1-\rho$ provides $\mu^{(\mathrm{Hill})}/\mu^{(\mathrm{Var})}\leq 1$. 
Equality holds if and only if 
$b_j/c_j=K$, a constant independent of $j$, {\it i.e.}~$b_j/c_j=b_1/c_1=1$. In other words, $| \mu^{(\mathrm{Hill})} | \leq | \mu^{(\mathrm{Var})} |$ 
with equality $\mu^{(\mathrm{Hill})} = \mu^{(\mathrm{Var})}$ if and only if either $\rho=0$ or $k_1/n_1=(k_j/n_j)(1+\operatorname{o}(1))$ for any $j$, meaning that the sample fraction in each machine should be (asymptotically) the same for asymptotic bias equality to hold.
\end{Rmk}
Corollary~\ref{coro:poolopt_mfixed_kinfty_bal} and Remark~\ref{rmk:bias} motivate the following result. 
\begin{Theo}
\label{theo:poolopt_mfixed_kinfty_bal2}
Under the conditions of Corollary~\ref{coro:poolopt_mfixed_kinfty_bal} with $\lambda\neq 0$, 
$
\sqrt{k} (\widehat{\gamma}_n(\widetilde{\boldsymbol{\omega}}_n^{(\mathrm{Var})})-\gamma) \stackrel{d}{\longrightarrow} \mathcal{N}( \lambda/(1-\rho), \gamma^2 )
$
if and only if $\rho=0$ or $k_j/n_j = (k/n)(1+\operatorname{o}(1))$ for any $j\in\{1,\ldots,m\}$. If moreover $k_j/n_j = (k/n)(1+\operatorname{O}(1/\sqrt{k}))$ for any $j$, 
we have in fact
$
\sqrt{k}(\widehat{\gamma}_n(\widetilde{\boldsymbol{\omega}}_n^{(\mathrm{Var})}) - \widehat{\gamma}_n^{(\mathrm{Hill})}(k)) = \operatorname{o}_{\mathbb{P}}(1).
$
\end{Theo}
This result states that adjusting the effective sample sizes
$k_j$ such that $k_j/n_j$ is constant across machines produces a variance-optimal distributed estimator that is asymptotically equivalent to the unfeasible Hill estimator built from the combined subsamples. This is much stronger than only sharing the same asymptotic distribution, referred to as the \emph{oracle} property by~\cite{chelizho2021} for the naive distributed estimator in the case $m \to\infty$ (compare also with the results of Section~\ref{sec:distrinf:allkbounded}).
In the particular sector of insurance, according to Supplement~B in~\cite{chelizho2021}, insurance companies may be willing to communicate and might effectively use the same sample fraction. Yet,
it is unlikely that such an adjustment can be performed in 
other sectors of activity,
since each machine will pick its own $k_j$ following an appropriate selection rule based only on its subsample. 

This difficulty with variance-optimal pooling motivates the focus on AMSE-optimal pooling in the frequent case of unequal sample fractions $k_j/n_j$. The following result compares the resulting AMSE-optimal distributed estimator with the benchmark Hill estimator.
\begin{Theo}
\label{theo:poolopt_mfixed_kinfty_AMSE}
Under the conditions and notation of Corollary~\ref{coro:pool_mfixed_kinfty_balanced}, set
\begin{align*}
\mbox{AMSE}^{(\mathrm{Hill})} &= \frac{1}{k} \left( \frac{\lambda^2}{(1-\rho)^2} + \gamma^2 \right) \mbox{ and } \mbox{AMSE}(\boldsymbol{\omega}) = \frac{1}{k} ((\boldsymbol{\omega}^{\top} \boldsymbol{B}_{\boldsymbol{c}})^2+\boldsymbol{\omega}^{\top} \mathbf{V}_{\boldsymbol{c}} \boldsymbol{\omega}) \\
\mbox{with } \boldsymbol{B}_{\boldsymbol{c}} &= \frac{\lambda}{1-\rho} \left( d_1^{\rho},\ldots, d_m^{\rho} \right)^{\top} \mbox{ and } \mathbf{V}_{\boldsymbol{c}} = \gamma^2 \left( \sum_{i=1}^m c_i^{-1} \right) \operatorname{diag}(c_1,\ldots,c_m). 
\end{align*}
Assume that the $b_j/c_j$ are not all equal to 1. Then $\mbox{AMSE}(\boldsymbol{\omega}^{(\mathrm{AMSE})}) \geq \mbox{AMSE}^{(\mathrm{Hill})}$ if and only if $|\lambda|\leq \lambda_0$, with 
\[
\lambda_0 = \gamma(1-\rho) \sqrt{\frac{S_{\rho}^2 - S_0^2}{S_0 S_{2\rho} - S_{\rho}^2}}, \mbox{ and } S_{\alpha}=\sum_{j=1}^m d_j^{\alpha}/c_j. 
\]
\end{Theo}
In the general situation where the marginal Hill estimators have unequal sample fractions $k_j/n_j$, it is remarkable that the AMSE-optimal distributed estimator can actually have a smaller AMSE than the benchmark Hill estimator itself under the necessary and sufficient condition that the bias component $| \lambda |$ is sufficiently large, as will be illustrated in Section~\ref{sec:fin:simul:distrinf}. 
It should also be noted 
that Theorem~\ref{theo:poolopt_mfixed_kinfty_AMSE} does not violate the minimax optimality property of the (benchmark) Hill estimator proved in~\cite{dre1998}
since it only states that the AMSE-optimal pooled estimator performs better than the Hill estimator within a certain class of heavy-tailed distributions. 

To estimate the AMSE-optimal weights in this distributed inference context, we assume 
as in Section~\ref{sec:main:optimal_EVI} that $A(t)=\gamma \beta t^{\rho}$ and note that 
\[
\lambda = \lim_{n\to\infty} \sqrt{k} A(n/k) = \gamma \beta \lim_{n\to\infty} \sqrt{k} \left( \frac{n}{k} \right)^{\rho} \mbox{ and } d_j=\lim_{n\to\infty} \frac{n_j k}{n k_j}. 
\]
The estimators $\widehat\beta_j$ and $\widehat\rho_j$ in Section~\ref{sec:main:optimal_EVI} of the second-order parameters $\beta\equiv\beta_j$ and $\rho\equiv\rho_j$ are restricted to each machine $j$ separately. We improve on these marginal estimators by using here the pooled versions $\widehat{\beta}_n(\boldsymbol{\omega}) = \sum_{j=1}^m \omega_j \widehat{\beta}_j$ and $\widehat{\rho}_n(\boldsymbol{\omega}) = \sum_{j=1}^m \omega_j \widehat{\rho}_j$, which in turn leads to the bias component estimator
$
\widehat{\lambda} = \widehat{\lambda}(\widehat{\boldsymbol{\omega}}_{\gamma},\widehat{\boldsymbol{\omega}}_{\beta},\widehat{\boldsymbol{\omega}}_{\rho}) = \widehat{\gamma}_n(\widehat{\boldsymbol{\omega}}_{\gamma}) \widehat{\beta}_n(\widehat{\boldsymbol{\omega}}_{\beta}) \times \sqrt{k} ( n/k )^{\widehat{\rho}_n(\widehat{\boldsymbol{\omega}}_{\rho})} ,
$
where $\widehat{\boldsymbol{\omega}}_{\gamma}$, $\widehat{\boldsymbol{\omega}}_{\beta}$ and $\widehat{\boldsymbol{\omega}}_{\rho}$ are three possibly random sets of weights. An obvious choice is $\omega_j=1/m$ for all three estimators. A more refined choice is variance-optimal weights $k_j/k$ for $\widehat{\boldsymbol{\omega}}_{\gamma}$ and $n_j/n$ for both $\widehat{\boldsymbol{\omega}}_{\beta}$ and $\widehat{\boldsymbol{\omega}}_{\rho}$ (recall that $\widehat{\beta}_j$ and $\widehat{\rho}_j$ use almost all the available observations in each machine). Set then 
\begin{align*}
\widetilde{\boldsymbol{B}}_{\boldsymbol{c}} &= \sqrt{ \sum_{j=1}^m k_j } \frac{\widehat{\gamma}_n(\widehat{\boldsymbol{\omega}}_{\gamma}) \widehat{\beta}_n(\widehat{\boldsymbol{\omega}}_{\beta})}{1-\widehat{\rho}_n(\widehat{\boldsymbol{\omega}}_{\rho})} \left( \left( \frac{n_1}{k_1} \right)^{\widehat{\rho}_n(\widehat{\boldsymbol{\omega}}_{\rho})}, \ldots, \left( \frac{n_m}{k_m} \right)^{\widehat{\rho}_n(\widehat{\boldsymbol{\omega}}_{\rho})} \right)^{\top} \\
\mbox{and } \widetilde{\mathbf{V}}_{\boldsymbol{c}} &=  \left( \sum_{j=1}^m k_j \right) \widehat{\gamma}_n^2(\widehat{\boldsymbol{\omega}}_{\gamma}) \operatorname{diag}(1/k_1,\ldots,1/k_m) ,
\end{align*}
and define $\widetilde{\boldsymbol{\omega}}_n^{(\mathrm{AMSE})}$ by replacing $\boldsymbol{B}_{\boldsymbol{c}}$ and $\mathbf{V}_{\boldsymbol{c}}$ in $\boldsymbol{\omega}^{(\mathrm{AMSE})}$ with 
$\widetilde{\boldsymbol{B}}_{\boldsymbol{c}}$ and $\widetilde{\mathbf{V}}_{\boldsymbol{c}}$. 
We 
obtain the following result as an immediate consequence of Corollary~\ref{coro:pool_mfixed_kinfty_balanced}. 
\begin{Coro}
\label{coro:pool_mfixed_kinfty_balanced_AMSE}
Work under the conditions of Corollary~\ref{coro:pool_mfixed_kinfty_balanced} with $\rho<0$ and $A(t)=\gamma \beta t^{\rho}$. Assume that, for all $j\in\{1,\ldots,m\}$, $\widehat{\beta}_j$ is a consistent estimator of $\beta$ and $(\widehat{\rho}_j-\rho) \log n_j=\operatorname{o}_{\mathbb{P}}(1)$. Then 
\[
\sqrt{k} (\widehat{\gamma}_n(\widetilde{\boldsymbol{\omega}}_n^{(\mathrm{AMSE})})-\gamma) \stackrel{d}{\longrightarrow} \mathcal{N}\left( \frac{\lambda}{1-\rho} \sum_{j=1}^m d_j^{\rho} \omega_j^{(\mathrm{AMSE})}, \gamma^2 \sum_{j=1}^m \frac{1}{c_j} \sum_{j=1}^m c_j (\omega_j^{(\mathrm{AMSE})})^2 \right).
\]
\end{Coro}
A confidence interval based on this AMSE-optimal estimator (or on the variance-optimal estimator) can then be constructed exactly as in Section~\ref{sec:main:inference}, and so we omit the details for the sake of brevity. This construction requires the knowledge of 
$k_j$, $n_j$, $\widehat{\beta}_j$, $\widehat{\rho}_j$ that each individual machine has to communicate to the central machine. 

We conclude this section 
by making use of the distributed inference-specific bias component estimator
to define the following bias-reduced versions for both the variance- and AMSE-optimal pooled estimators and derive their asymptotic distributions: 
\begin{align*}
\overline{\gamma}_n(\widetilde{\boldsymbol{\omega}}_n^{(\mathrm{Var})}) &= \widehat{\gamma}_n(\widetilde{\boldsymbol{\omega}}_n^{(\mathrm{Var})}) - \frac{1}{\sqrt{k}}(\widetilde{\boldsymbol{\omega}}_n^{(\mathrm{Var})})^{\top} \widetilde{\boldsymbol{B}}_{\boldsymbol{c}} \\
\mbox{and } \overline{\gamma}_n(\widetilde{\boldsymbol{\omega}}_n^{(\mathrm{AMSE})}) &= \widehat{\gamma}_n(\widetilde{\boldsymbol{\omega}}_n^{(\mathrm{AMSE})}) - \frac{1}{\sqrt{k}}(\widetilde{\boldsymbol{\omega}}_n^{(\mathrm{AMSE})})^{\top} \widetilde{\boldsymbol{B}}_{\boldsymbol{c}}.
\end{align*}
\begin{Coro}
\label{coro:pool_mfixed_kinfty_balanced_BR}
Under the conditions of Corollary~\ref{coro:pool_mfixed_kinfty_balanced_AMSE},
$\sqrt{k} (\overline{\gamma}_n(\widetilde{\boldsymbol{\omega}}_n^{(\mathrm{Var})})-\gamma) \stackrel{d}{\longrightarrow} \mathcal{N}(0, \gamma^2)$ and
\[
\sqrt{k} (\overline{\gamma}_n(\widetilde{\boldsymbol{\omega}}_n^{(\mathrm{AMSE})})-\gamma) \stackrel{d}{\longrightarrow} \mathcal{N}\left( 0, \gamma^2 \sum_{j=1}^m \frac{1}{c_j} \sum_{j=1}^m c_j (\omega_j^{(\mathrm{AMSE})})^2 \right).
\]
\end{Coro}

\subsection{Extreme quantile estimation}
\label{sec:distrinf:quantile}

We are now ready to compare the weighted geometric distributed estimator of an extreme quantile $q(1-p)$, as well as its variance- and AMSE-optimal versions, to the classical unfeasible Weissman estimator obtained directly from the combined subsamples, each defined as
\[
\widehat{q}_n^{\star}(1-p|\boldsymbol{\omega}) = \prod_{j=1}^m \left[ \left( \frac{k_j}{n_j p} \right)^{\widehat{\gamma}_j(k_j)} X_{n_j-k_j:n_j,j} \right]^{\omega_j}, \ \widehat{q}_n^{\star,(\mathrm{Hill})}(1-p|k) = \left( \frac{k}{n p} \right)^{\widehat{\gamma}_n^{(\mathrm{Hill})}(k)} X_{n-k:n}. 
\]
\begin{Coro}
\label{coro:poolopt_mfixed_kinfty_Weiss}
Work under the conditions of Corollary~\ref{coro:pool_mfixed_kinfty_balanced} with $\rho<0$. Pick $p=p(n)\to 0$ such that $k/(n p)\to \infty$ and $\sqrt{k}/\log(k/(n p))\to \infty$, as $n\to\infty$. Let $\boldsymbol{\omega}, \widehat{\boldsymbol{\omega}}_n$ be such that $\widehat{\boldsymbol{\omega}}_n^{\top} \boldsymbol{1}=1$ and $\widehat{\boldsymbol{\omega}}_n\stackrel{\mathbb{P}}{\longrightarrow} \boldsymbol{\omega}$. Then
\begin{align*}
    \frac{\sqrt{k}}{\log(k/(n p))} \left( \frac{\widehat{q}_n^{\star}(1-p|\widehat{\boldsymbol{\omega}}_n)}{q(1-p)} - 1 \right) &= \sqrt{k} (\widehat{\gamma}_n(\boldsymbol{\omega})-\gamma)+\operatorname{o}_{\mathbb{P}}(1) \\
    & \stackrel{d}{\longrightarrow} \mathcal{N}\left( \frac{\lambda}{1-\rho} \sum_{j=1}^m d_j^{\rho} \omega_j, \gamma^2 \sum_{j=1}^m \frac{1}{c_j} \sum_{j=1}^m c_j \omega_j^2 \right).
\end{align*}
If moreover $k_j/n_j = (k/n)(1+\operatorname{O}(1/\sqrt{k}))$ for any $j\in\{1,\ldots,m\}$, then $\widehat{q}_n^{\star}(1-p|\widetilde{\boldsymbol{\omega}}_n^{(\mathrm{Var})})$ is $\sqrt{k}/\log(k/(n p))-$asymptotically equivalent to $\widehat{q}_n^{\star,(\mathrm{Hill})}(1-p|k)$. 
Finally, if the $b_j/c_j$ are not all equal to 1, then under the conditions of Corollary~\ref{coro:pool_mfixed_kinfty_balanced_AMSE}, $\widehat{q}_n^{\star}(1-p|\widetilde{\boldsymbol{\omega}}_n^{(\mathrm{AMSE})})$ has a smaller AMSE than $\widehat{q}_n^{\star,(\mathrm{Hill})}(1-p|k)$ if and only if $|\lambda|>\lambda_0$, with the notation of Theorem~\ref{theo:poolopt_mfixed_kinfty_AMSE}.
\end{Coro}

\subsection{Extension to the case of at least one, but not all, very low $k_j$}
\label{sec:distrinf:onekbounded}

It may happen that the ratio $\max_{1\leq j\leq m} k_j/\min_{1\leq j\leq m} k_j$ 
is quite large, owing to uncertainty in data-driven selection rules. One may then want to simply discard the marginal estimates with a very low $k_j$ from the pooling procedure, but these estimates will have very low bias if all machines have comparable sample sizes, and hence it is more sensible to incorporate them into the distributed estimators of the tail index and extreme quantiles. From an asymptotic point of view, we obtain the following result for the variance-optimal distributed estimators in this situation of extremely unbalanced effective sample sizes. 
\begin{Theo}
\label{theo:pool_mfixed_kbounded}
Assume that condition $\mathcal{C}_2(\gamma,\rho,A)$ holds. Suppose that there is $\ell \in \{ 1,\ldots,m-1\}$ such that
on the one hand, for any $j\in \{1,\ldots,\ell\}$, $k_j\to\infty$ with $k_j/n_j\to 0$, $n_1/n_j\to b_j\in (0,\infty)$ and $k_1/k_j\to c_j\in (0,\infty)$ (with then $b_1,c_1=1$), and $\sqrt{k_j} A(n_j/k_j)\to \lambda_j \in \mathbb{R}$ as $n\to\infty$; and on the other hand, for any $j\in \{\ell+1,\ldots,m\}$, $k_j=k_j(n)$ is a nondecreasing sequence with $k_j/n_j\to 0$, $k_1/k_j\to \infty$ as $n\to\infty$ and $\sqrt{k_j} A(n_j/k_j) = \operatorname{O}(1)$.
Then,
if $\lambda=\lim_{n\to\infty} \sqrt{k} A(\sum_{i=1}^\ell n_i/\sum_{i=1}^{\ell} k_i) \in \mathbb{R}$ and $d_j=(c_j/b_j) \times (\sum_{i=1}^{\ell} c_i^{-1})/(\sum_{i=1}^{\ell} b_i^{-1})$, 
\[
\sqrt{k} (\widehat{\gamma}_n(\widetilde{\boldsymbol{\omega}}_n^{(\mathrm{Var})})-\gamma) \stackrel{d}{\longrightarrow} \mathcal{N}\left( \frac{\lambda}{1-\rho} \left( \sum_{j=1}^{\ell} \frac{1}{c_j} \right)^{-1} \left( \sum_{j=1}^{\ell} \frac{d_j^{\rho}}{c_j} \right), \gamma^2 \right).
\]
If moreover $\rho<0$ and $p=p(n)\to 0$ is such that $k/(n p)\to \infty$ and $\sqrt{k}/\log(k/(n p))\to \infty$, 
\begin{equation}
\label{eqn:condkbounded}
\max_{1\leq j\leq m} \frac{\sqrt{k_j}}{\sqrt{k}} \left| \frac{\log(k_j/(n_j p))}{\log(k/(np))} - 1 \right| \to 0 \mbox{ and } \max_{\ell+1\leq j\leq m} \frac{\sqrt{k_j}}{\sqrt{k}} \times \sqrt{k_j} |A(1/p)| = \operatorname{O}(1)
\end{equation}
as $n\to\infty$, then 
\[
\frac{\sqrt{k}}{\log(k/(n p))} \left( \frac{\widehat{q}_n^{\star}(1-p|\widetilde{\boldsymbol{\omega}}_n^{(\mathrm{Var})})}{q(1-p)} - 1 \right) \stackrel{d}{\longrightarrow} \mathcal{N}\left( \frac{\lambda}{1-\rho} \left( \sum_{j=1}^{\ell} \frac{1}{c_j} \right)^{-1} \left( \sum_{j=1}^{\ell} \frac{d_j^{\rho}}{c_j} \right), \gamma^2 \right).
\]
\end{Theo}
The variance-optimal distributed estimators 
therefore possess, under suitable conditions, the same asymptotic properties as if they were calculated without incorporating the machines using a very low $k_j$ into the pooling procedure. Note that the first part of condition~\eqref{eqn:condkbounded} automatically holds if in addition $n_1/n_j\to b_j\in (0,\infty)$ for any $j\in \{\ell+1,\ldots,m\}$, and the second part of the condition on $A(1/p)$ is obviously satisfied if, for any $j\in \{\ell+1,\ldots,m\}$, $k_j$ is bounded. 
Bias reduction can be similarly carried out by focusing only on the machines having a large $k_j$, so we omit the details for the sake of brevity. 

\subsection{The case of a large number of machines}
\label{sec:distrinf:allkbounded}

Our results above do not consider the case when all the $k_j$ are low, possibly even bounded, whose treatment is substantially different. When all the $k_j=k_j(n)$ are bounded in $n$, consistency of the distributed estimators necessarily requires a growing number of machines with $n$, namely $m=m(n)\to\infty$. Otherwise, when $m$ is fixed, the pooled estimators will contain only a bounded number of summands. In this context, we require the following fundamental assumption. 
\vskip1ex
\noindent
$(\mathcal{A})$ $m=m(n)\to\infty$ and the $n_j=n_j(n)$ satisfy $\inf_{1\leq j\leq m} n_j/\log m\to \infty$ as $n\to\infty$.
\vskip1ex
\noindent
This condition means that the amount of data stored in each machine grows with $n$, although the number of machines itself may grow at a much faster rate than the $n_j$. While the condition $m\to\infty$ ensures consistency, the condition $\inf_{1\leq j\leq m} n_j/\log m\to \infty$ is required to establish a precise control of the statistical errors 
arising in each machine. The proof is in this regard fundamentally different from the proofs of our theorems in the case of a fixed $m$. Another important difference is that the weight vector $\boldsymbol{\omega}\in \mathbb{R}^m$ now implicitly varies with $n$, so restrictions are needed to define a class of admissible weights. For example, the pooled estimator corresponding to the weight $(1,0,0,\ldots)$ is simply $\widehat{\gamma}_1(k_1)$, which is not consistent when $k_1$ is bounded. We thus introduce a {\it balanced allocation condition} on $\boldsymbol{\omega}=\boldsymbol{\omega}(n)$. 
\vskip1ex
\noindent
$(\mathcal{W})$ The weight vector $\boldsymbol{\omega}=\boldsymbol{\omega}(n)\in \mathbb{R}^m$ satisfies $\sum_{j=1}^m \omega_j=1$ as well as 
\[
\limsup_{n\to\infty} k \sum_{j=1}^m \frac{\omega_j^2}{k_j} <\infty \mbox{ and } \exists \delta>0, \ \frac{\sum_{j=1}^m k_j^{-\delta/2} (\omega_j^2/k_j)^{1+\delta/2}}{( \sum_{j=1}^m \omega_j^2/k_j )^{1+\delta/2}} \to 0 \mbox{ as } n\to\infty.
\]
This condition forbids $\boldsymbol{\omega}$ from having constant components (with respect to $n$), meaning that the weight should be roughly evenly spread out across machines. It also prevents $|\omega_j|$ from taking very large values despite summing up to 1, which corresponds to a stabilization condition. Any weight vector satisfying $m(n) \sup_{1\leq j\leq m} |\omega_j(n)|\leq \omega_0<\infty$, for all $n$, will automatically fulfill~$(\mathcal{W})$.
This encompasses the naive pooled estimator of $\gamma$, studied in~\cite{chelizho2021}. 
The 
distributed estimator $\widehat{\gamma}_n(\boldsymbol{\omega})$ has the following asymptotic properties.
\begin{Theo}
\label{theo:pool_minfty}
Under the conditions $(\mathcal{A})$, $(\mathcal{W})$ and $\mathcal{C}_2(\gamma,\rho,A)$, if $\sup_{1\leq j\leq m} k_j/n_j\to 0$ as $n\to\infty$ and $\sqrt{k} \sup_{1\leq j\leq m} |A(n_j/k_j)| = \operatorname{O}(1)$, then 
\[
\left( \sum_{j=1}^m \frac{\omega_j^2}{k_j} \right)^{-1/2} \left( \widehat{\gamma}_n(\boldsymbol{\omega})-\gamma-\frac{1}{1-\rho}\sum_{j=1}^m \omega_j k_j^{\rho} \frac{\Gamma(k_j-\rho+1)}{k_j!} A(n_j/k_j) \right) \stackrel{d}{\longrightarrow} \mathcal{N}(0,\gamma^2).
\]
Here $\Gamma$ denotes Euler's Gamma function. If also $k \sum_{j=1}^m \omega_j^2/k_j\to v\in [1,\infty)$ as $n\to\infty$, then 
\[
\sqrt{k} \left( \widehat{\gamma}_n(\boldsymbol{\omega})-\gamma-\frac{1}{1-\rho}\sum_{j=1}^m \omega_j k_j^{\rho} \frac{\Gamma(k_j-\rho+1)}{k_j!} A(n_j/k_j) \right) \stackrel{d}{\longrightarrow} \mathcal{N}(0,v\gamma^2).
\]
[Note that necessarily $k \sum_{j=1}^m \omega_j^2/k_j\geq 1$ by the Cauchy-Schwarz inequality.]
\end{Theo}
\begin{Rmk}[Comparison with previous results for fixed $m$]
\label{rmk:comparison} 
It is interesting to note that the asymptotic variance in Theorem~\ref{theo:pool_minfty} in fact matches the asymptotic variance in the 
case of fixed~$m$ in Corollary~\ref{coro:pool_mfixed_kinfty_balanced}, where the asymptotic variance is $\gamma^2$ multiplied by
\[
\sum_{j=1}^m \frac{1}{c_j} \sum_{j=1}^m c_j \omega_j^2 = \lim_{n\to\infty}  \sum_{j=1}^m \frac{k_j}{k_1} \sum_{j=1}^m \frac{k_1}{k_j}\omega_j^2 = \lim_{n\to\infty} k \sum_{j=1}^m \frac{\omega_j^2}{k_j} 
\]
because $c_j=\lim_{n\to\infty} k_1/k_j$. In the case $m\to\infty$, the right-hand side is nothing but $v$. By contrast, the asymptotic bias term is substantially different. 
\end{Rmk}
\begin{Rmk}[Comparison with~\cite{chelizho2021}]
\label{rmk:comparisonchen} 
Theorem~\ref{theo:pool_minfty} revisits and generalizes Theorems~1,~2 and~3 of~\cite{chelizho2021} in different directions: we deal with generic weighted distributed estimation instead of a naive pooling, we unravel not only the case of unbounded but also bounded (possibly unbalanced) $k_j$ in a single unified result, we handle the more realistic case of unbalanced sample sizes $n_j$, all of this under the natural weaker version $(\mathcal{A})$ of Condition~A in~\cite{chelizho2021}. We also remark that the asymptotic variance $\gamma^2$ obtained in Theorem~2 of~\cite{chelizho2021}, with the naive weights $\omega_j=1/m$ for $j\in \{ 1,\ldots,m\}$, is not correct. We fix the problem by proving that the asymptotic variance is in fact $v \gamma^2$, where in general $v>1$. This higher variance is not surprising because, in the case of unbalanced $k_j$, machines with the lowest $k_j$ tend to provide less information than those with the largest $k_j$, and therefore a loss of information should be expected in comparison with the case where all the $k_j$ are equal.
This insight can be checked by considering, for example, a simple situation where $X$ is purely Pareto distributed with tail index $\gamma$, the number $m$ of machines is even, and $k_j=1$ for $j$ odd and $k_j=2$ for $j$ even. In this situation, each $\widehat{\gamma}_j(k_j)$ is in fact simply a sum of $k_j$ independent exponential random variables with mean $\gamma$ and variance $\gamma^2$, and hence 
\[
\operatorname{Var}\left( \sum_{j=1}^m \sqrt{k} (\widehat{\gamma}_n(1/m,\ldots,1/m) - \gamma) \right) = \frac{3m/2}{m^2} \left( \sum_{l=1}^{m/2} \gamma^2 + \sum_{l=1}^{m/2} \frac{\gamma^2}{2} \right) 
= \frac{9\gamma^2}{8}.
\]
This matches our result since $\sum_{j=1}^m k_j=3m/2$ and $\sum_{j=1}^m \omega_j^2/k_j=3/4m$, so that $v=9/8$. 
As such, the distributed Hill estimator with $\omega_j=1/m$ typically does not achieve the so-called ``oracle property'' claimed in~\cite{chelizho2021} even if $\sqrt{k} \sup_{1\leq j\leq m} |A(n_j/k_j)| \to 0$ or $\rho=0$. In their terminology, this property holds when the distributed Hill estimator inherits the same speed of convergence and asymptotic distribution as the infeasible benchmark Hill estimator.
\end{Rmk}
Remark~\ref{rmk:comparisonchen} motivates the following corollary on the estimator with $\omega_j=\widehat{\omega}_j^{(\mathrm{Var})}=k_j/k$, which clearly satisfies condition~$(\mathcal{W})$.
\begin{Coro}
\label{coro:pool_mfixed_allkbounded}
Under conditions $(\mathcal{A})$ and $\mathcal{C}_2(\gamma,\rho,A)$, if $\sup_{1\leq j\leq m} k_j/n_j\to 0$ as $n\to\infty$ and $\sqrt{k} \sup_{1\leq j\leq m} |A(n_j/k_j)| = \operatorname{O}(1)$, then 
\[
\sqrt{k} \left( \widehat{\gamma}_n(\widehat{\boldsymbol{\omega}}_n^{(\mathrm{Var})})-\gamma-\frac{1}{1-\rho}\sum_{j=1}^m \frac{k_j}{k} \times k_j^{\rho} \frac{\Gamma(k_j-\rho+1)}{k_j!} A(n_j/k_j) \right) \stackrel{d}{\longrightarrow} \mathcal{N}(0,\gamma^2).
\]
\end{Coro}
Therefore, as in the case of bounded $m$, the distributed estimator with $\boldsymbol{\omega}=\widehat{\boldsymbol{\omega}}_n^{(\mathrm{Var})}$ is asymptotically variance-optimal. It is this weighted estimator which possesses the ``oracle property''.

We now turn to the asymptotic behavior of the geometrically weighted extreme quantile estimator $\widehat{q}_n^{\star}(1-p|\boldsymbol{\omega})$. 
Perhaps surprisingly, this distributed estimator is in fact generally inconsistent when the $k_j$ are bounded, {\it i.e.}~$\limsup_{n\to\infty} \sup_{1\leq j\leq m} k_j(n) <\infty$. The rationale is that, while the bias of the individual shape parameter estimators $\widehat{\gamma}_j(k_j)$ is small, and so averaging them out as $\widehat{\gamma}_n(\boldsymbol{\omega})$ creates a consistent estimator as the number of machines increases, the individual scale parameter estimators $X_{n_j-k_j:n_j,j}$ are fundamentally biased estimators of $q(1-k_j/n_j)$ when $k_j$ is 
fixed. 
As such, the Weissman extrapolation of $X_{n_j-k_j:n_j,j}$ to the far tail in conjunction with $\widehat{\gamma}_n(\boldsymbol{\omega})$ can no longer be correctly applied. The workaround to ensure consistency of $X_{n_j-k_j:n_j,j}$ is to choose $k_j\to\infty$ with $k_j/n_j\to 0$ as $n\to\infty$. With these growing $k_j$, the distributed extreme quantile estimator is, as expected, consistent and asymptotically normal. These insights are summarized in the following result. 
\begin{Theo}
\label{theo:pool_minfty_Weiss}
Work under the conditions of Theorem~\ref{theo:pool_minfty} with $\rho<0$. Pick $p=p(n)\to 0$ such that $k/(n p)\to \infty$ and $\sqrt{k}/\log(k/(n p))\to \infty$ as $n\to\infty$, and assume that 
\[
\sup_{1\leq j\leq m} \left| \frac{\log(k_j/(n_j p))}{\log(k/(np))} - 1 \right| \to 0 \mbox{ as } n\to\infty.
\]
\begin{enumerate}[label=(\roman*), wide, labelindent=0pt]
\item If the $k_j=k_j(n)$ are bounded, i.e.~$\limsup_{n\to\infty} \sup_{1\leq j\leq m} k_j(n) <\infty$, and $\omega_j\geq 0$ for any $j$, then 
$
\widehat{q}_n^{\star}(1-p|\boldsymbol{\omega})/q(1-p) 
$
does not converge to 1 in probability.
\item If the $k_j=k_j(n)$ are such that $\inf_{1\leq j\leq m} k_j\to \infty$ and $k \sum_{j=1}^m \omega_j^2/k_j\to v\in [1,\infty)$ as $n\to\infty$, then
\begin{multline*}
    \frac{\sqrt{k}}{\log(k/(n p))} \left( \frac{\widehat{q}_n^{\star}(1-p|\boldsymbol{\omega})}{q(1-p)} - 1 - \frac{1}{1-\rho}\sum_{j=1}^m \omega_j A(n_j/k_j) \right) \\
    = \sqrt{k} \left( \widehat{\gamma}_n(\boldsymbol{\omega})-\gamma-\frac{1}{1-\rho}\sum_{j=1}^m \omega_j A(n_j/k_j) \right)+\operatorname{o}_{\mathbb{P}}(1) \stackrel{d}{\longrightarrow} \mathcal{N}(0,v\gamma^2).
\end{multline*}
\end{enumerate}
\end{Theo}
Given that the $k_j$ are uniformly bounded, the assumption $k/(np)\to\infty$ is exactly $m/(np)\to\infty$, which is guaranteed to hold in the typical case $np\to c\in [0,\infty)$ of interest in extreme value analysis. The extra assumption on the $k_j$ and $n_j$ compared to Theorem~\ref{theo:pool_minfty} ensures that the effective sample fractions $k_j/n_j$ are not too dissimilar across machines; it is satisfied if, for instance, $0<\liminf_{n\to\infty} (n/k) \inf_{1\leq j\leq m} k_j/n_j \leq \limsup_{n\to\infty} (n/k) \sup_{1\leq j\leq m} k_j/n_j<\infty$. A weaker version of this condition already appears in Theorem~\ref{theo:pool_mfixed_kbounded}. 

\section{Filtering to handle dependence and covariates}
\label{sec:filtering}

In practice, the data are often recorded with relevant covariates, or are stationary but weakly dependent in a way that can be modeled by a standard time series. Besides, when the $X_j$ share the same tail index, they typically have also asymptotically proportional extreme quantiles (see Remark~\ref{rmk:hetero}). This suggests that $(X_1,\ldots,X_m)$ can be 
modeled in many situations by a general location-scale model 
\begin{equation}
\label{eqn:locscale}
X_j = g_j(\bm{Z}_j) + \sigma_j(\bm{Z}_j) \varepsilon_j, \ 1\leq j\leq m
\end{equation}
where the unobserved noise vector $\boldsymbol{\varepsilon}=(\varepsilon_1,\ldots,\varepsilon_m)$ has marginal tail quantile functions $U_j$ satisfying the conditions $\mathcal{C}_2(\gamma \boldsymbol{1},\boldsymbol{\rho},\boldsymbol{A})$ and $(\mathcal{H})$, and its bivariate survival copulae $\overline{C}_{j,\ell}$ satisfy $\mathcal{J}(\boldsymbol{R})$. The functions $g_j$ and $\sigma_j>0$ are unknown measurable functions of $\bm{Z}_j\in \mathbb{R}^{l_j}$, for some $l_j\geq 1$. The covariates $\bm{Z}_j$ can be fully observed (in traditional regression settings) or partially or not at all observed (in a time series model which includes past unobserved innovations or volatility terms). The noise variable $\varepsilon_j$ is assumed to be independent of $\bm{Z}_j$.
\vskip1ex
\noindent
Let then the pairs $(X_{i,j}, \bm{Z}_{i,j})$, $1\leq i\leq n_j$, be part of a strictly stationary sequence such that
$X_{i,j} = g_j(\bm{Z}_{i,j}) + \sigma_j(\bm{Z}_{i,j}) \varepsilon_{i,j}$ for $1\leq j\leq m$ and $1\leq i\leq n_j$. The $\boldsymbol{\varepsilon}_i=(\varepsilon_{i,1},\ldots,\varepsilon_{i,m})$ are assumed to be independent copies of $\boldsymbol{\varepsilon}$ as above. A reasonable idea to eliminate the heteroskedasticity and dependence in the data $(X_{i,j})$ is to first estimate the location and scale components $g_j$ and $\sigma_j$ of the model (under suitable identifiability and regularity conditions), and then filter the data to obtain residuals $\widehat{\varepsilon}_{i,j}^{(n_j)}$ close to the unobserved errors $\varepsilon_{i,j}$. This results in $j$ residual-based Hill estimators
$
\widehat{\gamma}_j(k_j) = k_j^{-1} \sum_{i=1}^{k_j} \log(\widehat{\varepsilon}_{n_j-i+1:n_j,j}^{(n_j)}/\widehat{\varepsilon}_{n_j-k_j:n_j,j}^{(n_j)}).
$
These can be combined in a pooled version $\widehat{\gamma}_n(\boldsymbol{\omega})=\sum_{j=1}^m \omega_j \widehat{\gamma}_j(k_j)$ whose asymptotic normality can be proved under a high-level condition on the discrepancy between $\widehat{\varepsilon}_{i,j}^{(n_j)}$ and $\varepsilon_{i,j}$.
\begin{Theo}
\label{theo:pool_mfixed_kinfty_filtering}
Assume that $\boldsymbol{\varepsilon}=(\varepsilon_1,\ldots,\varepsilon_m)$ satisfies assumptions $\mathcal{C}_2(\gamma \boldsymbol{1},\boldsymbol{\rho},\boldsymbol{A})$ and $\mathcal{J}(\boldsymbol{R})$. Under the conditions of Theorem~\ref{theo:pool_mfixed_kinfty_general} on $k_j$, $n_j$ and $\boldsymbol{\omega}$, if 
\[
\max_{1\leq j \leq m} \sqrt{k_j} \max_{1\leq i\leq n_j} \frac{|\widehat{\varepsilon}_{i,j}^{(n_j)} - \varepsilon_{i,j}|}{1+|\varepsilon_{i,j}|} \stackrel{\mathbb{P}}{\longrightarrow} 0 ,
\]
then
$
\sqrt{k} (\widehat{\gamma}_n(\boldsymbol{\omega})-\gamma) \stackrel{d}{\longrightarrow} \mathcal{N}( \boldsymbol{\omega}^{\top} \boldsymbol{B}_{\boldsymbol{c}}, \boldsymbol{\omega}^{\top} \mathbf{V}_{\boldsymbol{c}} \boldsymbol{\omega} )
$
with $\boldsymbol{B}_{\boldsymbol{c}}$ and $\mathbf{V}_{\boldsymbol{c}}$ defined analogously as in Theorem~\ref{theo:pool_mfixed_kinfty_general}.
If $\widehat{\boldsymbol{\omega}}_n^{\top} \boldsymbol{1}=1$ with $\widehat{\boldsymbol{\omega}}_n\stackrel{\mathbb{P}}{\longrightarrow} \boldsymbol{\omega}$, then $\sqrt{k}(\widehat{\gamma}_n(\widehat{\boldsymbol{\omega}}_n) - \widehat{\gamma}_n(\boldsymbol{\omega}))=\operatorname{o}_{\mathbb{P}}(1)$. In particular, $\widehat{\gamma}_n(\widehat{\boldsymbol{\omega}}_n)$ has the same $\sqrt{k}-$asymptotic behavior as $\widehat{\gamma}_n(\boldsymbol{\omega})$.
\end{Theo}
\begin{Rmk}[On the importance of filtering without pooling residuals]
\label{rmk:regression}
One might argue that pooling directly the residuals themselves in a sample and applying the traditional Hill estimator may be more efficient. However, if the model is misspecified, heteroskedasticity could still remain in the residuals, and those 
with the largest scale might swamp the other residuals in the pooled sample, resulting thus in a large loss of estimation accuracy. Pooling the residual-based Hill estimates instead provides more protection against departures from the assumed location-scale model. Related points in the standard (non-pooling) context are discussed in Remarks~2 and~3 in~\cite{girstuauc2021aos}.
\end{Rmk}
The condition on the discrepancy between $\widehat{\varepsilon}_{i,j}^{(n_j)}$ and $\varepsilon_{i,j}$ in Theorem~\ref{theo:pool_mfixed_kinfty_filtering} is typically satisfied as soon as the location and scale components $g_j$ and $\sigma_j$ are estimated at a faster rate than $\sqrt{k_j}$. This can easily be checked theoretically in a variety of regression models, see~\cite{girstuauc2021aos} for examples. The presence of the denominator $1+|\varepsilon_{i,j}|$ in the condition makes it also possible to handle heteroskedasticity. 
\vskip1ex
\noindent
The ultimate interest in conditional extreme value modeling under~\eqref{eqn:locscale} is to estimate extreme $(1-p)$th quantiles of $X_j$ given $\bm{Z}_j=\bm{z}_j$, defined as
$
q_{X_j|\bm{Z}_j=\bm{z}_j}(1-p) = g_j(\bm{z}_j) + \sigma_j(\bm{z}_j) q_j(1-p)
$
by location equivariance and positive homogeneity of quantiles, where $q_j$ is the quantile function of~$\varepsilon_j$. Given consistent estimators $\widehat{g}_j(\bm{z}_j)$ of $g_j(\bm{z}_j)$ and $\widehat{\sigma}_j(\bm{z}_j)$ of~$\sigma_j(\bm{z}_j)$, and using the Weissman estimator of $q_j(1-p)$ from the residuals $\widehat\varepsilon_{i,j}^{(n_j)}$, 
one can then estimate $q_{X_j|\bm{Z}_j=\bm{z}_j}(1-p)$ by
\[
\widehat{q}_{X_j|\bm{Z}_j=\bm{z}_j}^{\star}(1-p|k_j,\boldsymbol{\omega}) = \widehat{g}_j(\bm{z}_j) + \widehat{\sigma}_j(\bm{z}_j) \left( \frac{k_j}{n_j p} \right)^{\widehat{\gamma}_n(\boldsymbol{\omega})} \widehat{\varepsilon}_{n_j-k_j:n_j,j}^{(n_j)}. 
\]
When $\boldsymbol{\varepsilon}$ satisfies the condition $(\mathcal{H})$ described in Section~\ref{sec:main:extquant}, all associated quantiles $q_j(1-p)$ become asymptotically equivalent and can then be estimated by a geometrically pooled estimator, which leads to the following location-scale estimator of $q_{X_j|\bm{Z}_j=\bm{z}_j}(1-p)$:
\[
\widetilde{q}_{X_j|\bm{Z}_j=\bm{z}_j}^{\star}(1-p|\boldsymbol{\omega}) = \widehat{g}_j(\bm{z}_j) + \widehat{\sigma}_j(\bm{z}_j) \prod_{j=1}^m \left[ \left( \frac{k_j}{n_j p} \right)^{\widehat{\gamma}_j(k_j)} \widehat{\varepsilon}_{n_j-k_j:n_j,j}^{(n_j)} \right]^{\omega_j}.
\]
Our final asymptotic result establishes the asymptotic normality of these two estimators. 
\begin{Theo}
\label{theo:pool_mfixed_kinfty_filtering_Weiss}
Work under the conditions of Theorem~\ref{theo:pool_mfixed_kinfty_filtering} with $\rho_j<0$ for all $j\in \{1,\ldots,m\}$. Pick $p=p(n)\to 0$ such that $k/(n p)\to \infty$ and $\sqrt{k}/\log(k/(n p))\to \infty$ as $n\to\infty$. Let $\boldsymbol{\omega}, \widehat{\boldsymbol{\omega}}_n$ be such that $\widehat{\boldsymbol{\omega}}_n^{\top} \boldsymbol{1}=1$ and $\widehat{\boldsymbol{\omega}}_n\stackrel{\mathbb{P}}{\longrightarrow} \boldsymbol{\omega}$. Finally, assume that the estimators $\widehat{g}_j(\bm{z}_j)$ and $\widehat{\sigma}_j(\bm{z}_j)$ satisfy $\widehat{g}_j(\bm{z}_j) - g_j(\bm{z}_j) = \operatorname{O}_{\mathbb{P}}(1)$ and $\sqrt{k_j} (\widehat{\sigma}_j(\bm{z}_j) - \sigma_j(\bm{z}_j)) = \operatorname{O}_{\mathbb{P}}(1)$. Then, for any $j$, 
\[
\frac{\sqrt{k}}{\log(k_j/(n_j p))} \left( \frac{\widehat{q}_{X_j|\bm{Z}_j=\bm{z}_j}^{\star}(1-p|k_j,\widehat{\boldsymbol{\omega}}_n)}{q_{X_j|\bm{Z}_j=\bm{z}_j}(1-p)} - 1 \right) = \sqrt{k} (\widehat{\gamma}_n(\boldsymbol{\omega}) - \gamma) + \operatorname{o}_{\mathbb{P}}(1)
\] 
which converges weakly to $\mathcal{N}( \boldsymbol{\omega}^{\top} \boldsymbol{B}_{\boldsymbol{c}}, \boldsymbol{\omega}^{\top} \mathbf{V}_{\boldsymbol{c}} \boldsymbol{\omega} )$. If moreover $(\mathcal{H})$ holds then, for any $j$, 
\[
\frac{\sqrt{k}}{\log(k/(n p))} \left( \frac{\widetilde{q}_{X_j|\bm{Z}_j=\bm{z}_j}^{\star}(1-p|\widehat{\boldsymbol{\omega}}_n)}{q_{X_j|\bm{Z}_j=\bm{z}_j}(1-p)} - 1 \right) = \sqrt{k} (\widehat{\gamma}_n(\boldsymbol{\omega}) - \gamma) + \operatorname{o}_{\mathbb{P}}(1).
\]
\end{Theo}

\section{Finite-sample study}
\label{sec:fin}
\subsection{Simulation experiments}
\label{sec:fin:simul}

We investigate the finite-sample performance of our proposed inferential methodology, first in the general pooling framework for heavy-tailed distributions (Section~\ref{sec:main}) and then in the distributed inference framework (Section~\ref{sec:distrinf}). To save space we only report a brief description of our simulated models and conclusions.

\subsubsection{General setup: Pooling for tail index and extreme quantile inference}\label{sec:fin:simul:genpool}

Dimensions $m\in \{2,3,4,5\}$ were considered, with balanced, weakly unbalanced and strongly unbalanced sample sizes. Our statistical models 
had either unit Fr\'echet or absolute Student-$t$ ({\it i.e.}~the absolute value of a Student-$t$) marginal distributions with 1 degree of freedom. The dependence structure between 
margins was given by four copulae: the Clayton and Gumbel (Archimedean) copulae, and the Gaussian and Student copulae. The Clayton and Gaussian copulae are cases of asymptotic independence, while the Gumbel and Student copulae are cases of asymptotic dependence. 
All marginal distributions had equal tail indices $\gamma_j=1$.

Our first experiment, for a total sample size of $n=1{,}000$ across all subsamples, compares four pooled tail index estimators (naive, variance-optimal, AMSE-optimal, and AMSE-optimal with pooled second-order estimates as in the comment below Corollary~\ref{coro:pool_mfixed_kinfty_general_BR}) with the benchmark Hill estimator applied to the pooled dataset on each of the aforementioned models. We also compare the related four geometrically pooled extreme quantile estimators $\widehat{q}_n^{\star}(1-p|\widehat{\boldsymbol{\omega}}_n)$ at level $1-p=0.999$ with the naive arithmetic mean of the Weissman estimators $\widehat{q}_j^{\star}(1-p|k_j)$ in each subsample and the benchmark Weissman estimator applied to the pooled dataset. We compute Monte Carlo approximations of the Mean Squared Error (MSE) and of the actual coverage probability for the asymptotic confidence intervals with $95\%$ nominal level arising from our asymptotic theory, see Corollaries~\ref{coro:pool_testing} and~\ref{coro:pool_testing_quantile} for our proposed estimators; for the Hill (resp.~Weissman) estimator, we assume that the asymptotic distribution is normal with mean 0 and variance $\gamma^2/k$ (resp.~$\gamma^2\times \log^2(k/(np))$), see Theorem~3.2.5, p.74 in~\cite{haafer2006} (resp.~Theorem~4.3.8, p.138 therein). The variance-optimal and AMSE-optimal estimators outperform by far the naive pooling estimator on the basis of the MSE, when there is strong unbalance between sample sizes. Differences in performance get larger as the unbalance increases. They also perform comparably to the Hill estimator on pooled data. At the inferential level, confidence intervals deduced using the variance-optimal and AMSE-optimal estimators are typically substantially narrower than those provided using naive pooling. When there is asymptotic independence between samples, the pooling methods and the benchmarks on pooled data both provide confidence intervals having correct coverage. This is no longer the case when substantial dependence is present, with the $95\%$ confidence intervals constructed using the benchmark Hill estimator having actual coverage that can be as low as $75\%$. Conclusions about extreme quantile estimation are similar, with the added fact that geometrically pooled estimators outperform by far the naive arithmetic mean of Weissman estimators.

\subsubsection{Distributed inference of extreme values}
\label{sec:fin:simul:distrinf}

In our second experiment, we assume that the marginal distributions are i.i.d.~Fr\'echet, absolute Student-$t$ or Burr distributed, with tail index $\gamma=1$. We consider dimensions $m\in \{ 5,10,20 \}$ in balanced and highly unbalanced setups. We compared again the tail index and extreme quantile estimators described in Section~\ref{sec:fin:simul:genpool}. The results suggest that our proposed variance-optimal and AMSE-optimal methods (with pooled second-order estimates) perform comparably to the unfeasible Hill and Weissman estimators applied to the pooled dataset and outperform the naive distributed estimators, with shorter confidence intervals having correct coverage, and lower MSE. Geometric pooling is clearly beneficial as far as extreme quantile estimation is concerned. The gain of using our proposed distributed estimators increases as the unbalance between sample sizes increases. The AMSE-optimal distributed estimator is also overall the best when sample fractions are unequal and substantially different.

\subsection{Data analysis}
\label{sec:fin:data}

We discuss two concrete applications of our methodology to insurance and rainfall data.

\subsubsection{Distributed inference for car insurance data}
\label{sec:fin:data:distr}

\begin{figure}[t!]
\centering
\includegraphics[width=0.32\textwidth]{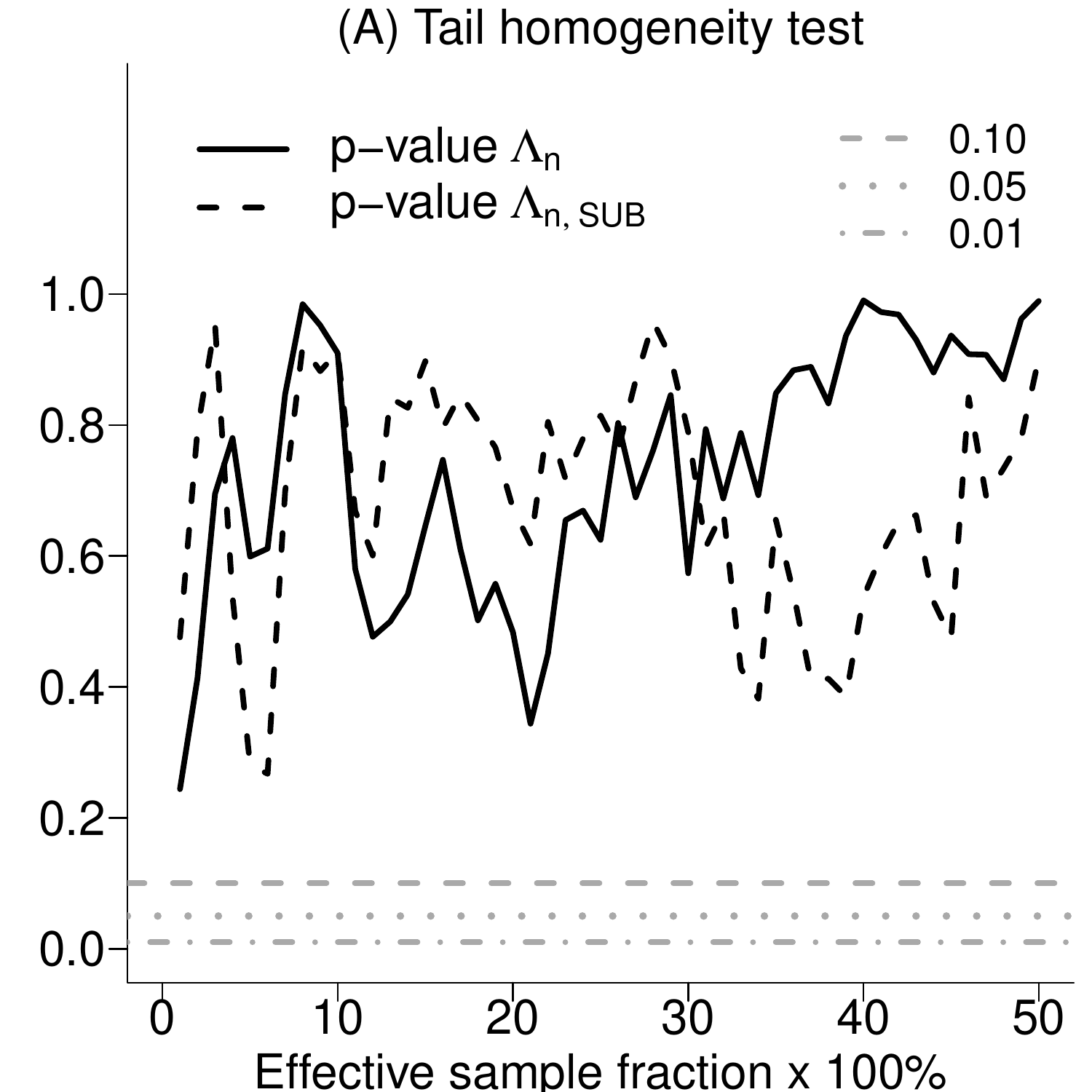}
\includegraphics[width=0.32\textwidth]{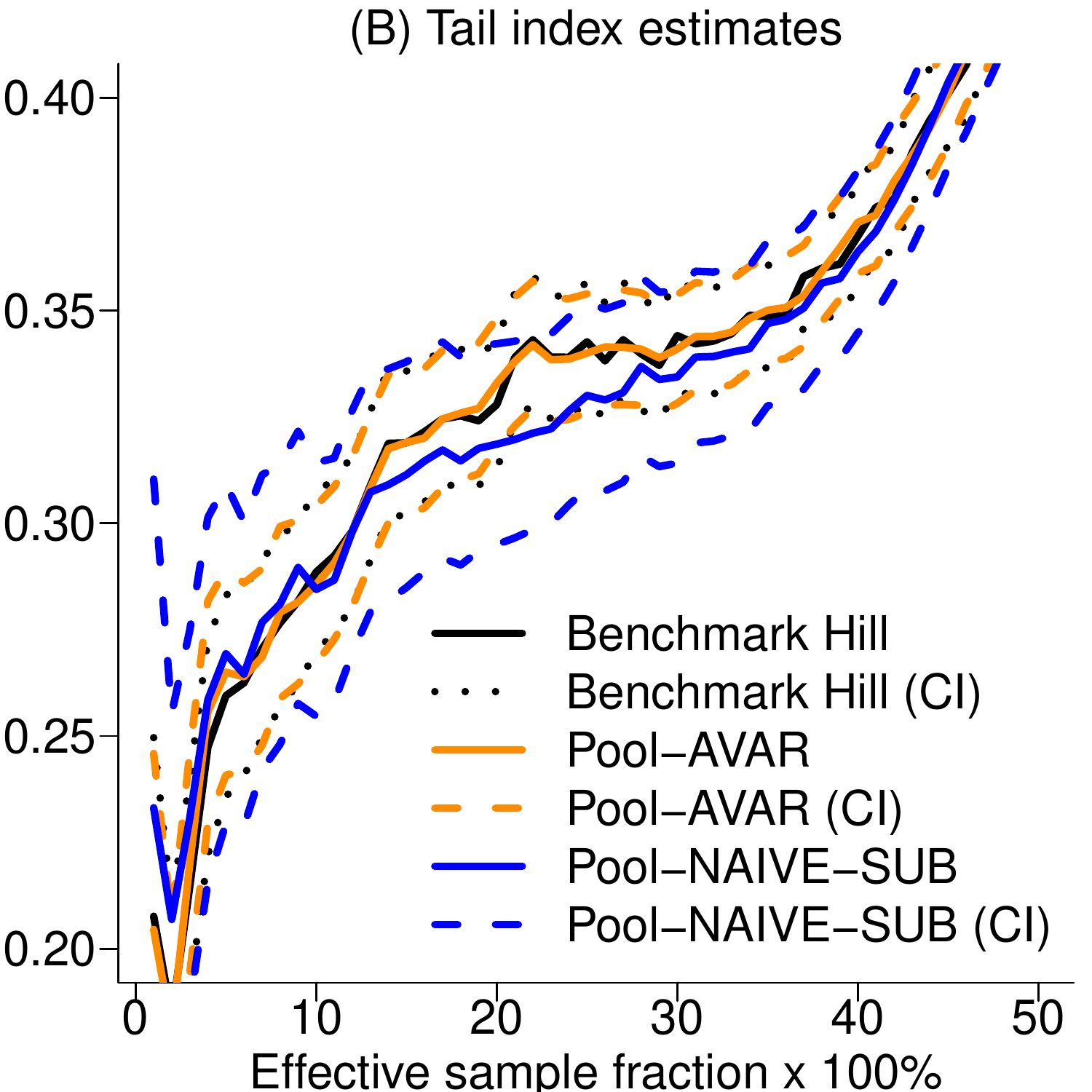}
\includegraphics[width=0.32\textwidth]{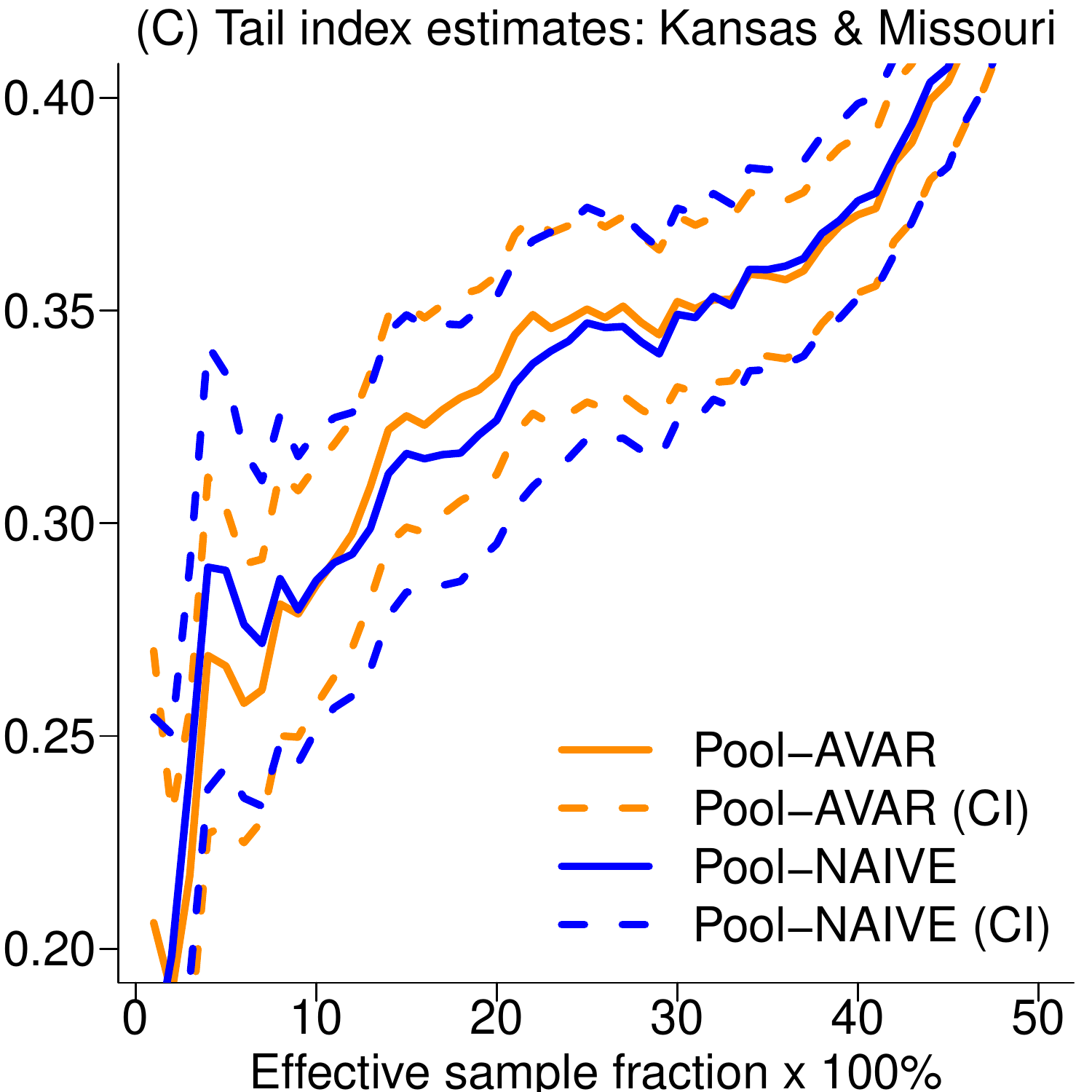} \\[10pt]
\includegraphics[width=0.32\textwidth]{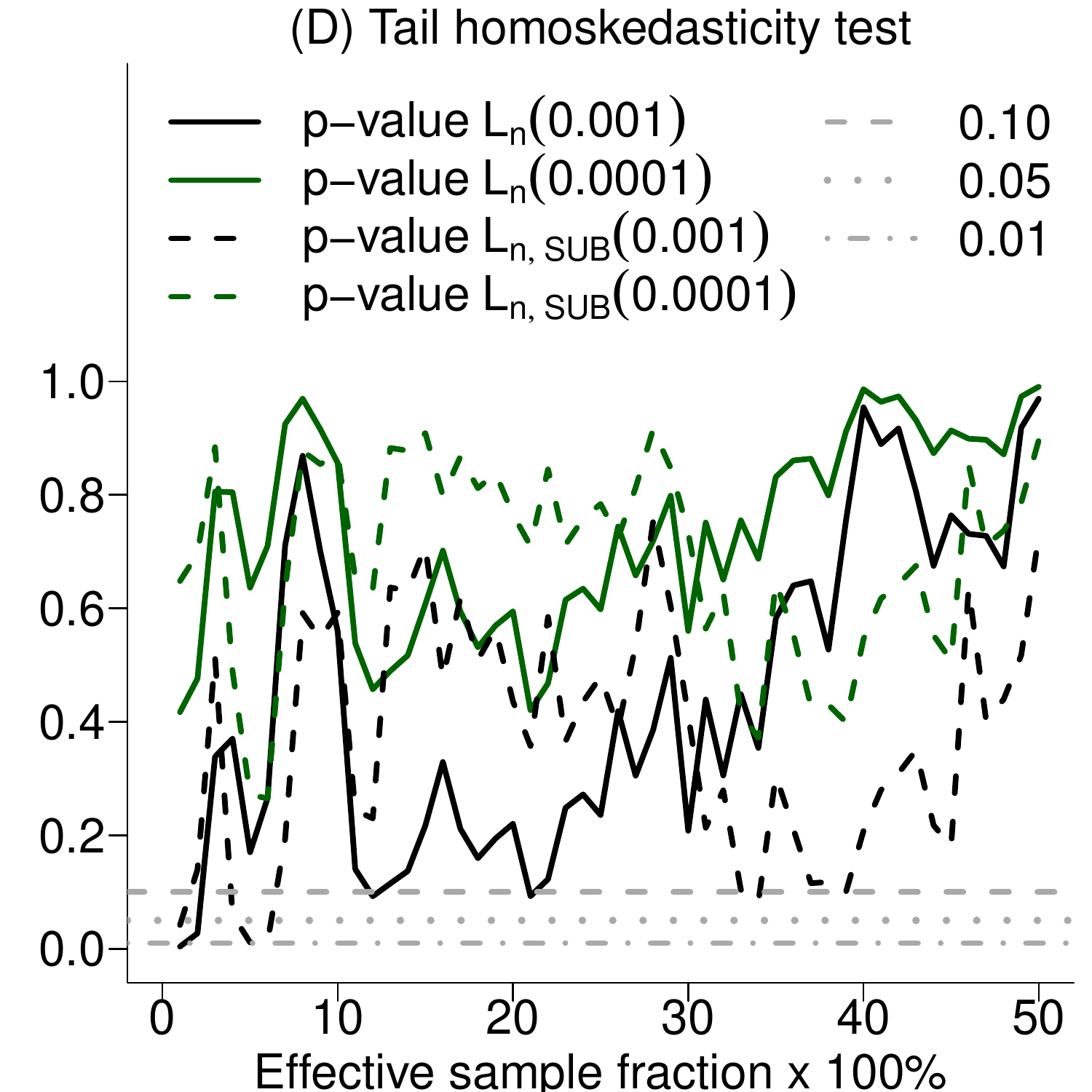}
\includegraphics[width=0.32\textwidth]{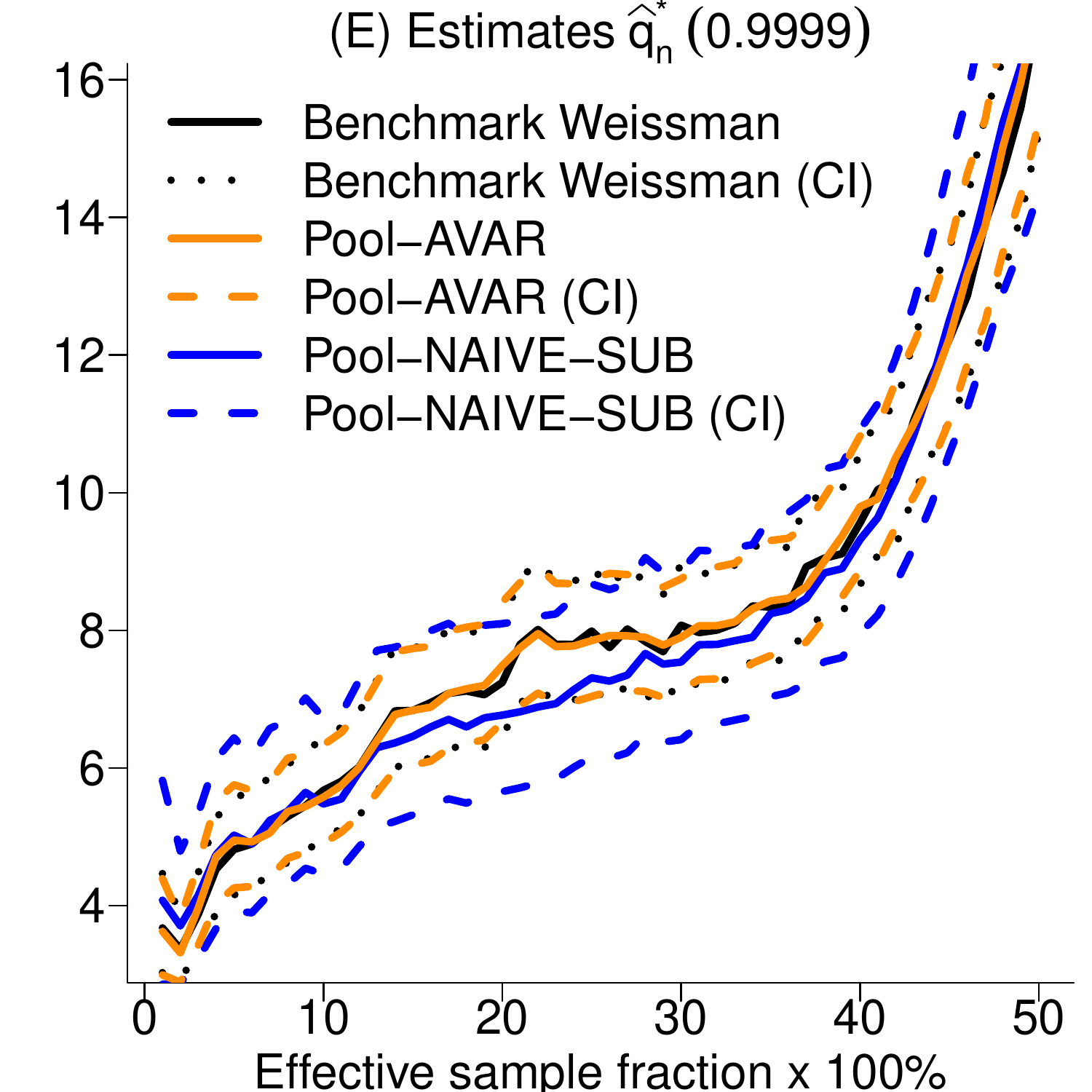}
\includegraphics[width=0.32\textwidth]{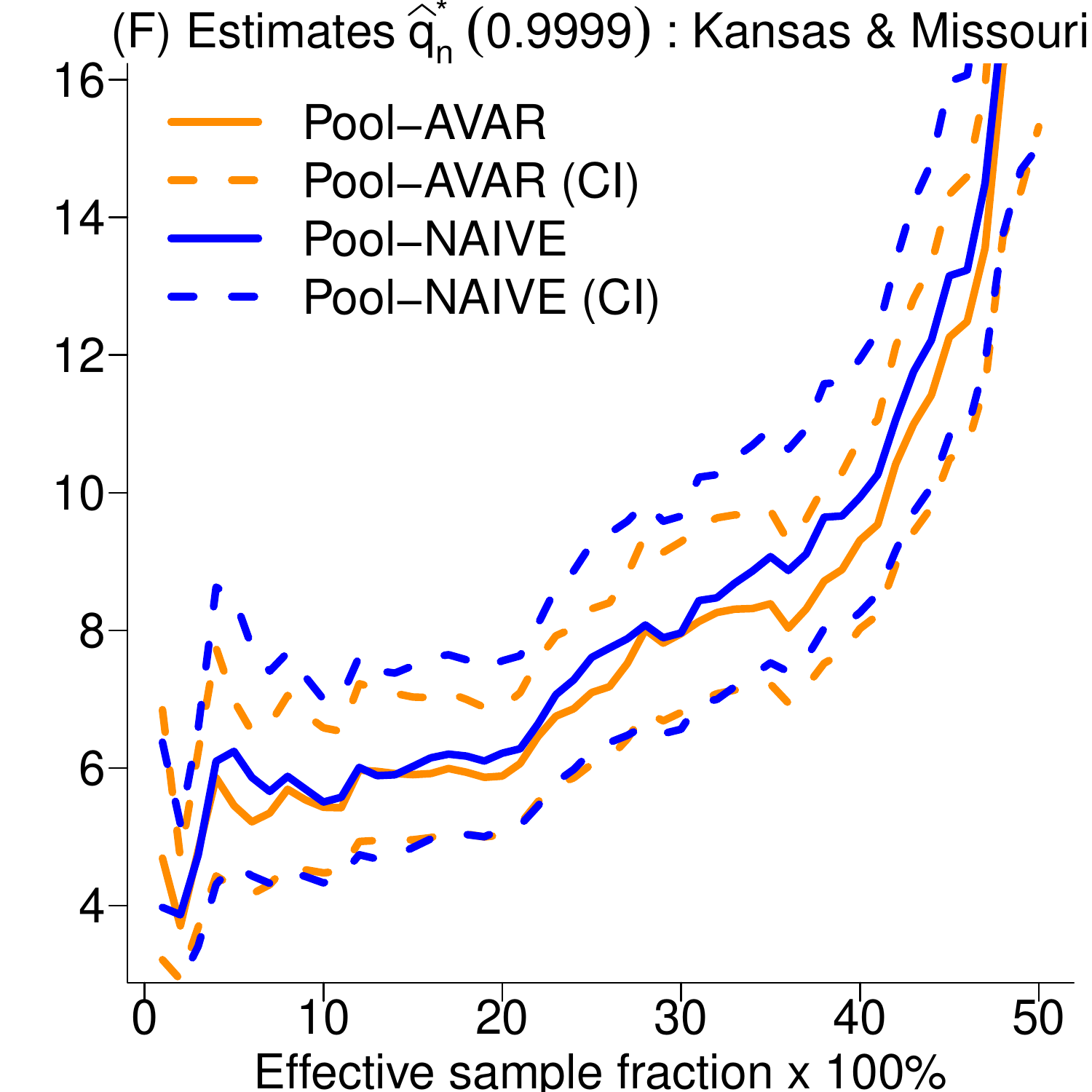}
\caption{Car insurance data. In (A) and (D), $\Lambda_{n,\textrm{SUB}}$ and $L_{n,\textrm{SUB}}(p)$ denote the test statistics $\Lambda_n$ and $L_n(p)$ calculated from the subsampled data of total size $5 \times 700 = 3{,}500$. In (B) and (C), Pool-AVAR and Pool-NAIVE respectively denote the variance-optimal pooled estimator and the naive pooled estimator of the tail index. In (E) and (F), Pool-AVAR and Pool-NAIVE respectively denote the variance-optimal pooled quantile estimator $\widehat{q}_n^{\star}(1-p|\widehat{\boldsymbol{\omega}}^{(\mathrm{Var})})$ and its unweighted analog. Pool-AVAR-SUB and Pool-NAIVE-SUB stand for these estimators calculated on the subsampled data. All estimates are represented as functions of the sample fraction $k_j/n_j$, assumed to be identical for each $j$.}
\label{fig:carpool:results}
\end{figure}

The first dataset comprises total claim amounts for car insurance companies in the five US states of Iowa ($n_1=2{,}601$), Kansas ($n_2=798$), Missouri ($n_3=3{,}150$), Nebraska ($n_4=1{,}703$), and Oklahoma ($n_5=882$) between January and February 2011, for a total sample size of $n=9{,}134$. 
The choice of this dataset follows the same setup as in~\cite{chelizho2021} who assume that each company
cannot share its data, making the calculation of Hill and Weissman estimates from pooled data inapplicable, but each company is willing to share its statistical analysis to enhance its appraisal of tail risk. Unlike~\cite{chelizho2021}, however, our distributed inference method can handle the different subsample sizes $n_j$ and hence the full dataset, and allows to estimate extreme quantiles. Their analysis requires first a subsampling step to guarantee the same subsample sizes. We compare our results using the full data with those obtained from their method after subsampling at random 700 observations in each state, as described in Supplement~B of~\cite{chelizho2021}. As a benchmark, we use the hypothetical Hill and Weissman estimates, obtained directly from the combined $n$ data points. Results are given in Figure~\ref{fig:carpool:results}.

We first check the equality of tail indices by testing for tail homogeneity across the 5 states on full data and the subsampled data in each state, using the theory developed in Section~\ref{sec:main:inference} under the constraint of independence between subsamples (see Remark~\ref{rmk:asyindep}). 
From the p-values corresponding to our test statistic $\Lambda_n$ in Figure~\ref{fig:carpool:results}(A), we can comfortably conclude the equality of individual tail indices at the three significance levels $0.10,\, 0.05$ and $0.01$. It is remarkable that the p-values plot remains quite stable when moving from the full 5 samples of total size $9{,}134$ to the subsamples of total size $5 \times 700 = 3{,}500$.
This indicates that the asymptotic chi-square 
regime is attained reasonably quickly. 

Figure~\ref{fig:carpool:results}(B) compares our variance-optimal distributed  
estimator $\widehat{\gamma}_n(\widehat{\boldsymbol{\omega}}_n^{(\mathrm{Var})})$, based on the full data, with the naive distributed estimator $\widehat\gamma_n(1/m,\ldots,1/m)$ of~\cite{chelizho2021}, which relies on the subsampled data, and with the benchmark Hill estimator $\widehat{\gamma}_n^{(\mathrm{Hill})}$ along with their respective asymptotic $95\%$ confidence intervals. 
In contrast to the naive estimates 
and their associated confidence intervals, 
our optimal weighted estimates 
and their confidence intervals 
are, respectively, 
almost identical to
the 
Hill estimates 
and their corresponding confidence intervals, 
as is to be expected from Theorem~\ref{theo:poolopt_mfixed_kinfty_bal2}. Our variance-optimal confidence intervals are found to be around $40\%$ shorter than those of~\cite{chelizho2021}. We arrive at a similar conclusion, in Figure~\ref{fig:carpool:results}(C), when restricting the analysis to the branches in Kansas and Missouri, whose subsample sizes $798$ and $3{,}150$ are strongly imbalanced, and using the full data from these states for both the variance-optimal and naive distributed estimates.
Here, the variance-optimal confidence intervals are found to be roughly $20\%$ shorter than those relative to the naive estimator.
The test theory of extreme quantile equivalence developed in Corollary~\ref{coro:pool_testing_quantile} is implemented for the two extreme quantile levels $1-p=0.999\approx 1-1/\max_j n_j$ and $1-p=0.9999\approx 1-1/n$, 
resulting in the p-values from the test statistic $L_n(p)$ displayed in Figure~\ref{fig:carpool:results}(D) for both full and subsampled data. 
The test 
overall allows to accept the assumption of tail homoskedasticity
across states, with p-values getting higher as $p$ decreases. The rationale behind this behavior in this distributed setting is that, as extreme quantile levels increase, the shape of the approximating Pareto distribution gets more important relative to its scale. As such, because mere differences in scale can no longer be detected in the far tail as $p\downarrow 0$, the test actually becomes less powerful against the sub-alternative of proportional quantiles. 
Finally, the resulting variance-optimal distributed estimates and confidence intervals for extreme quantiles are found to be virtually indistinguishable from the ideal Weissman analogs, whereas they appreciably outperform the naive distributed competitors, as can be seen 
from Figure~\ref{fig:carpool:results}(E) and~(F)
for $p=0.0001$. 

\subsubsection{Pooling for regional inference on extreme rainfall}
\label{sec:fin:data:rain}

\begin{figure}[t!]
\centering
        \includegraphics[width=0.5\textwidth]{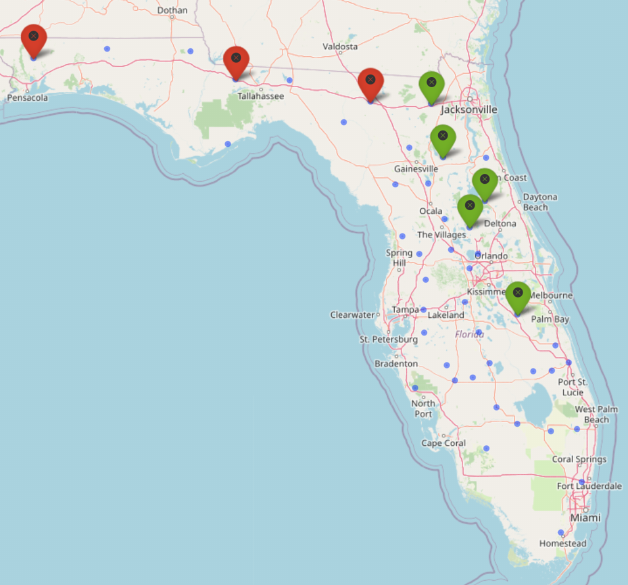} \\[10pt]
        \includegraphics[width=0.32\textwidth]{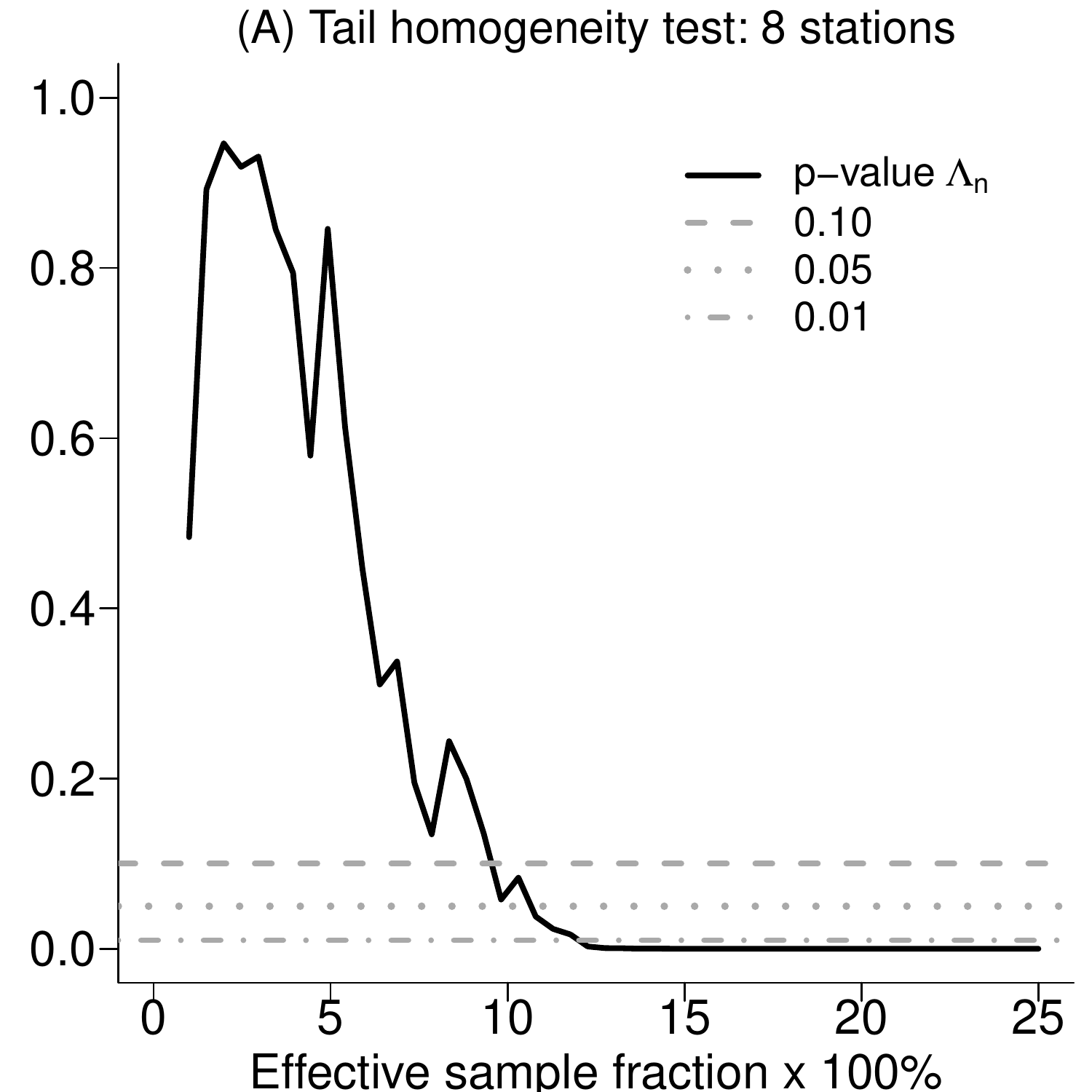} 
        \includegraphics[width=0.32\textwidth]{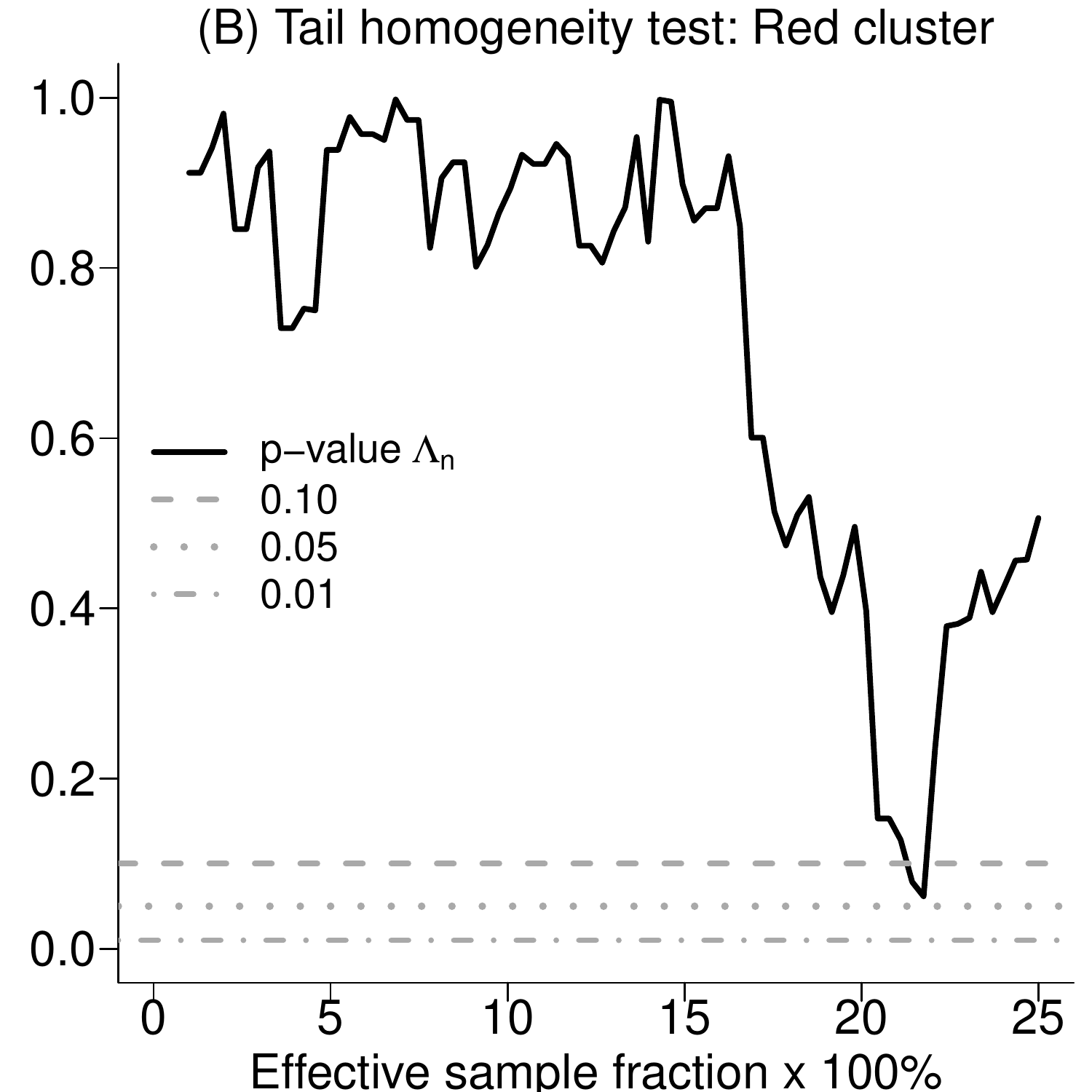} 
        \includegraphics[width=0.32\textwidth]{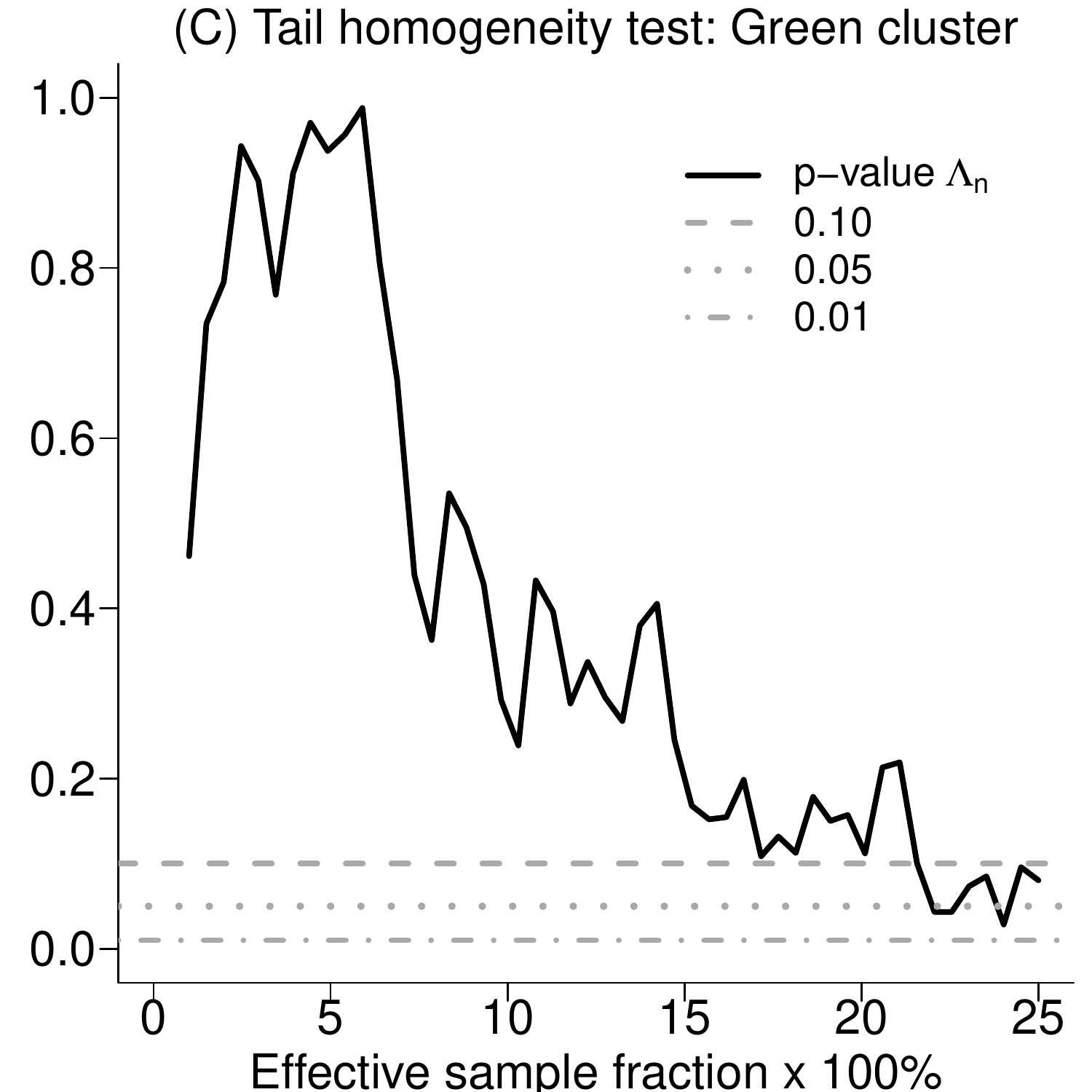} \\[10pt]
        \includegraphics[width=0.32\textwidth]{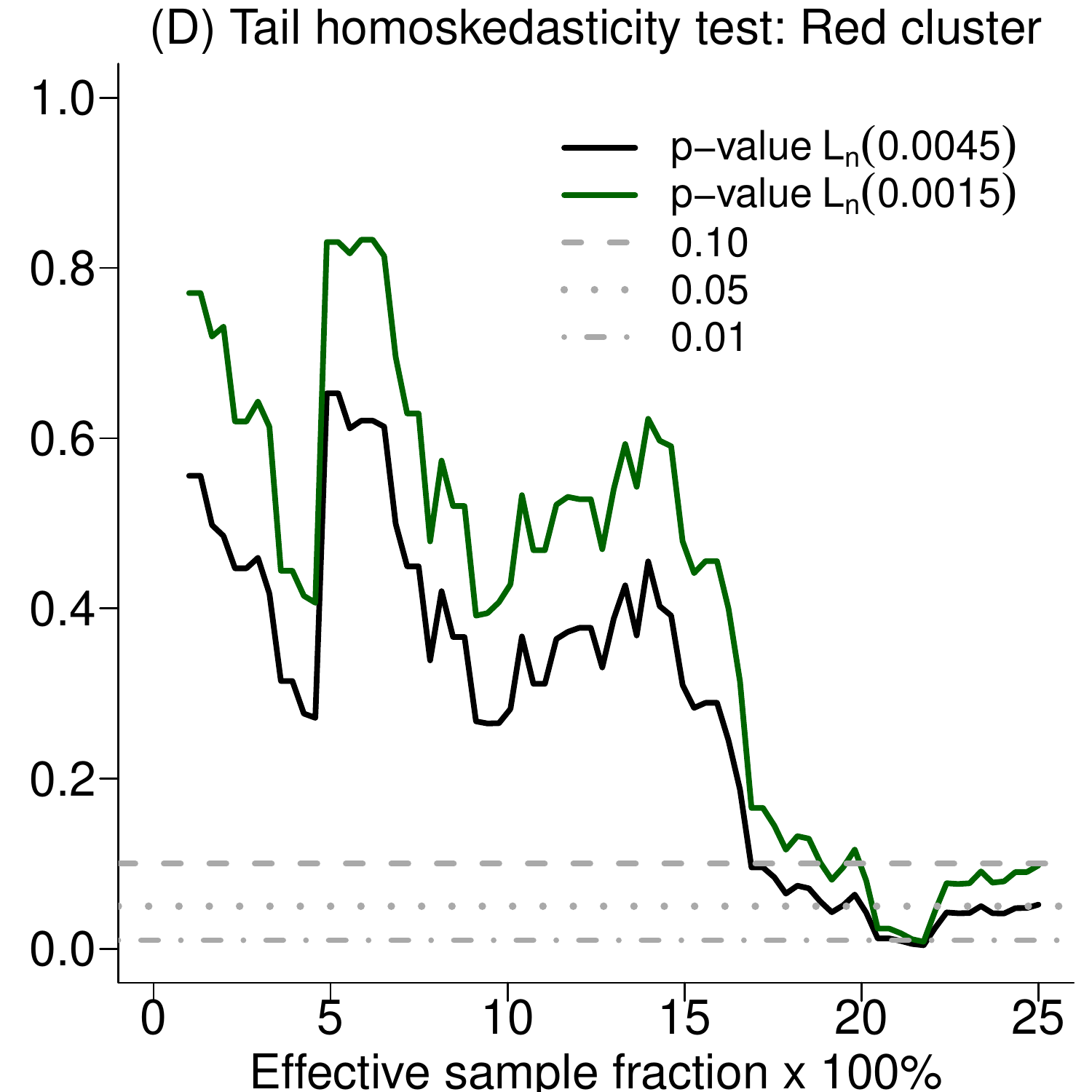} 
        \includegraphics[width=0.32\textwidth]{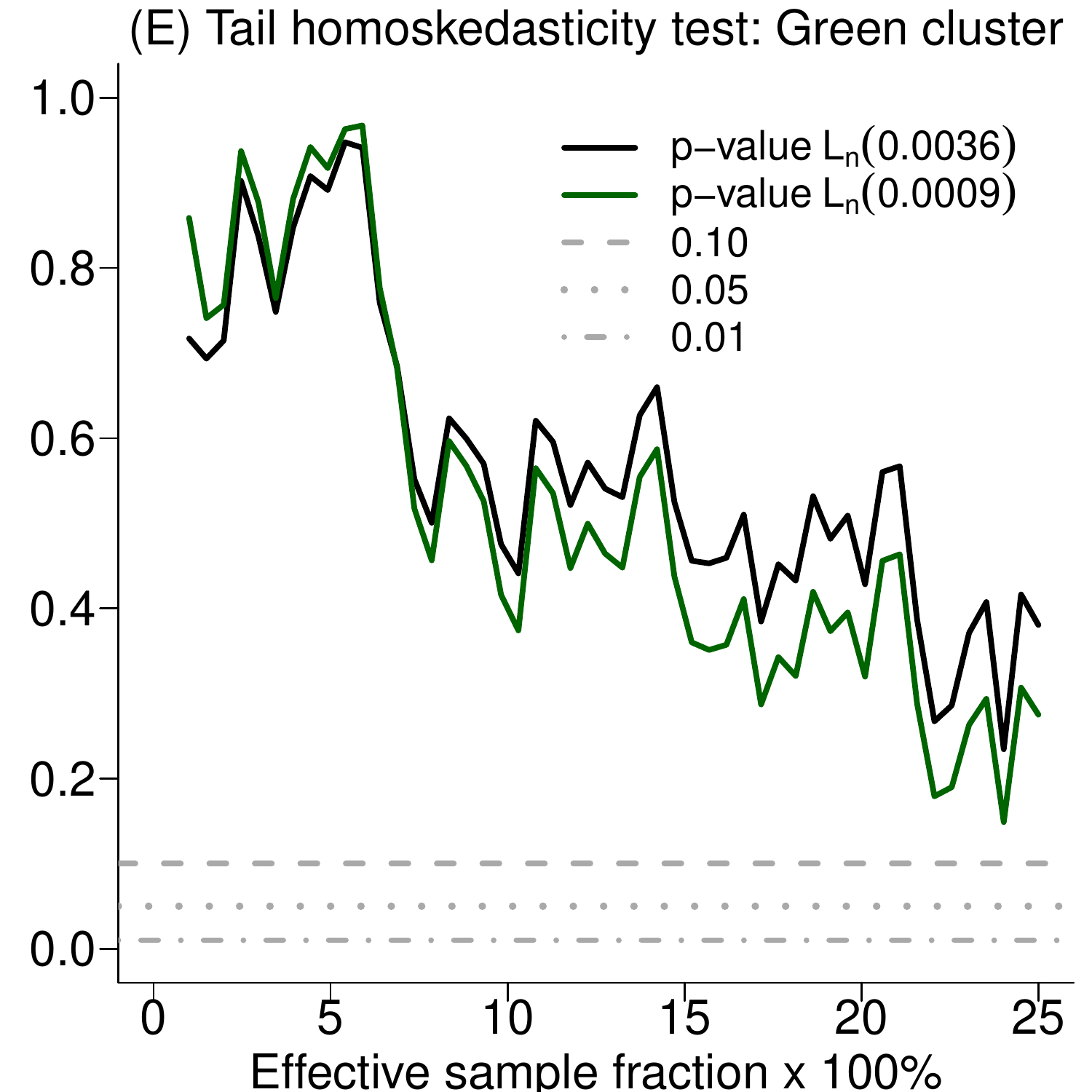} 
        \includegraphics[width=0.32\textwidth]{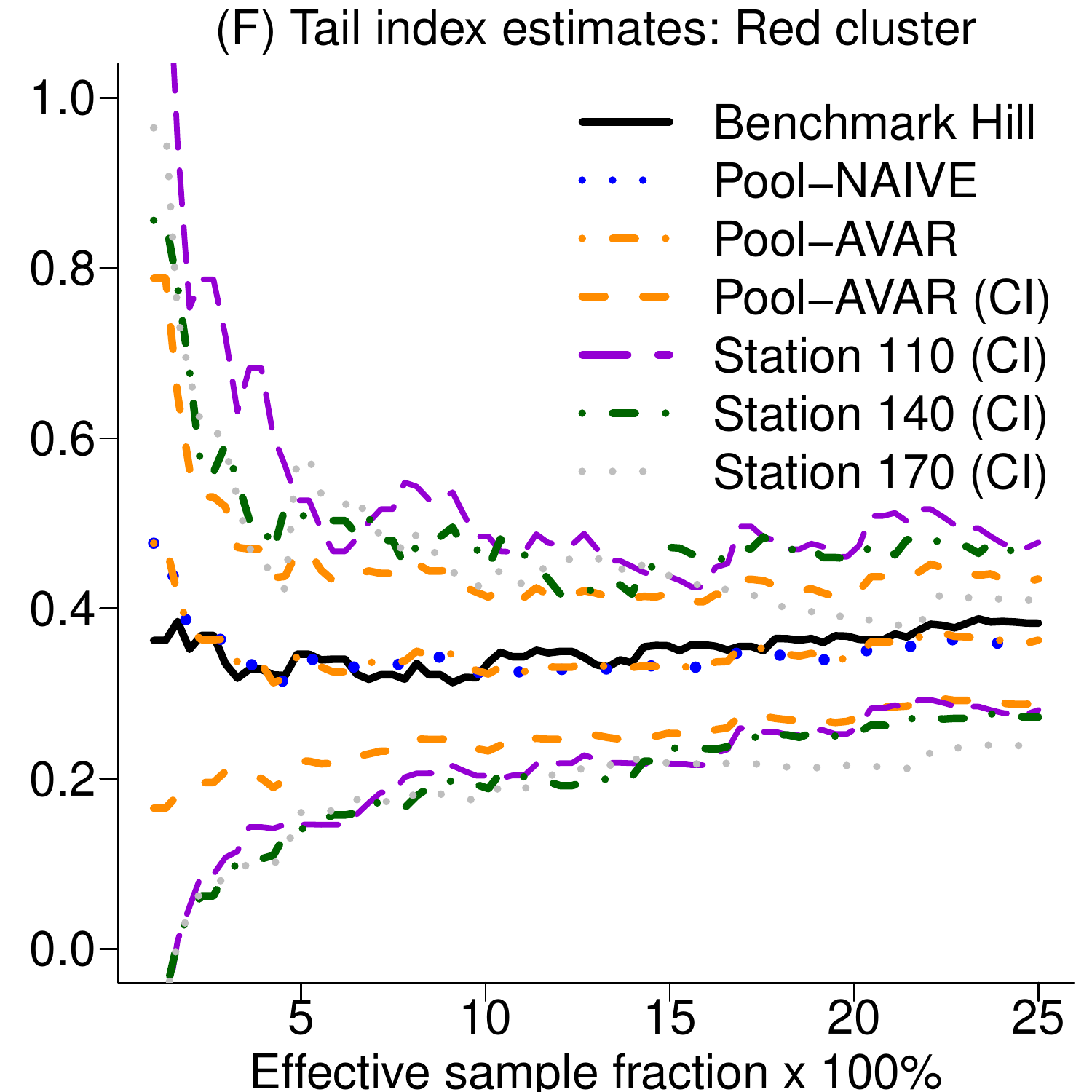} \\[10pt]
        \includegraphics[width=0.32\textwidth]{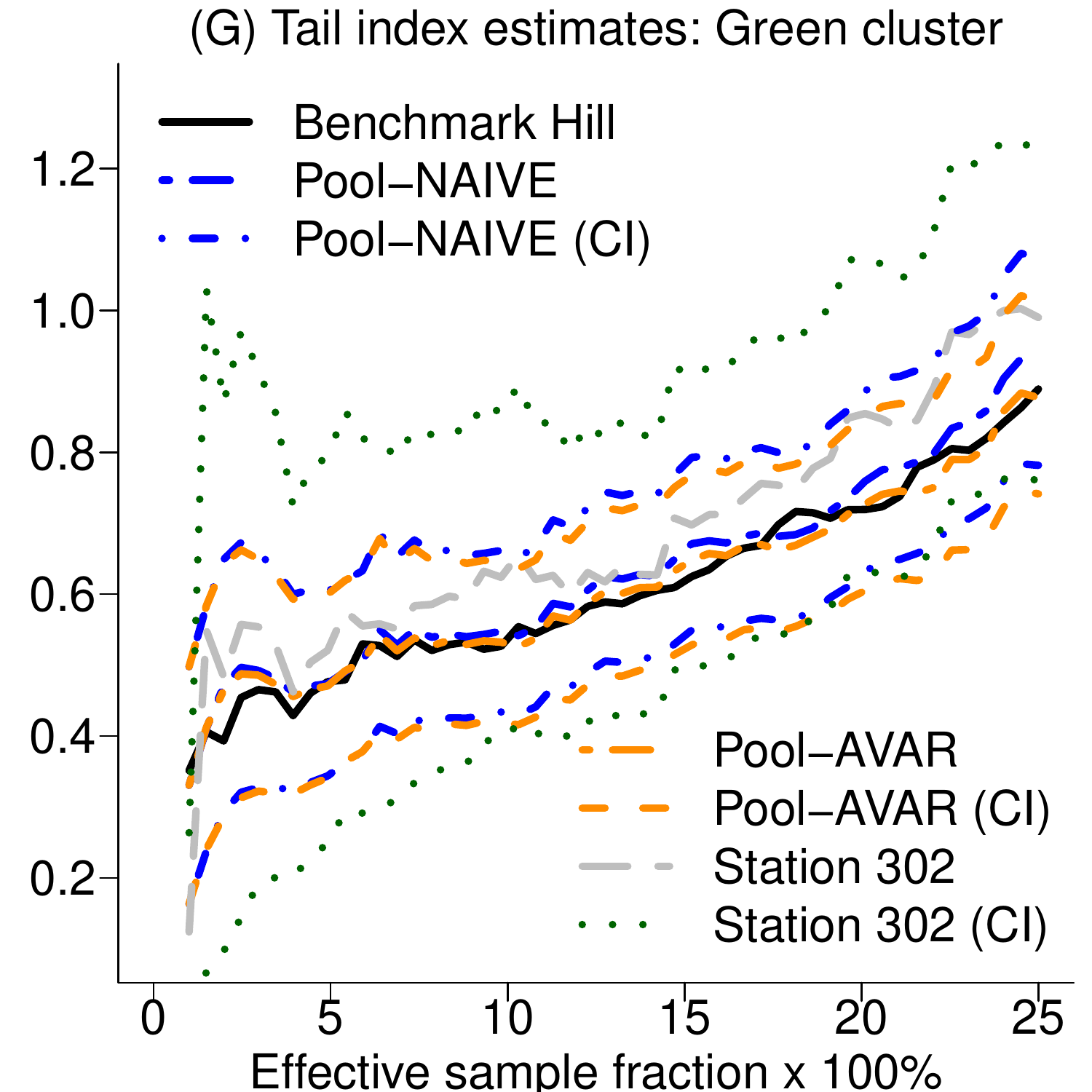} 
        \includegraphics[width=0.32\textwidth]{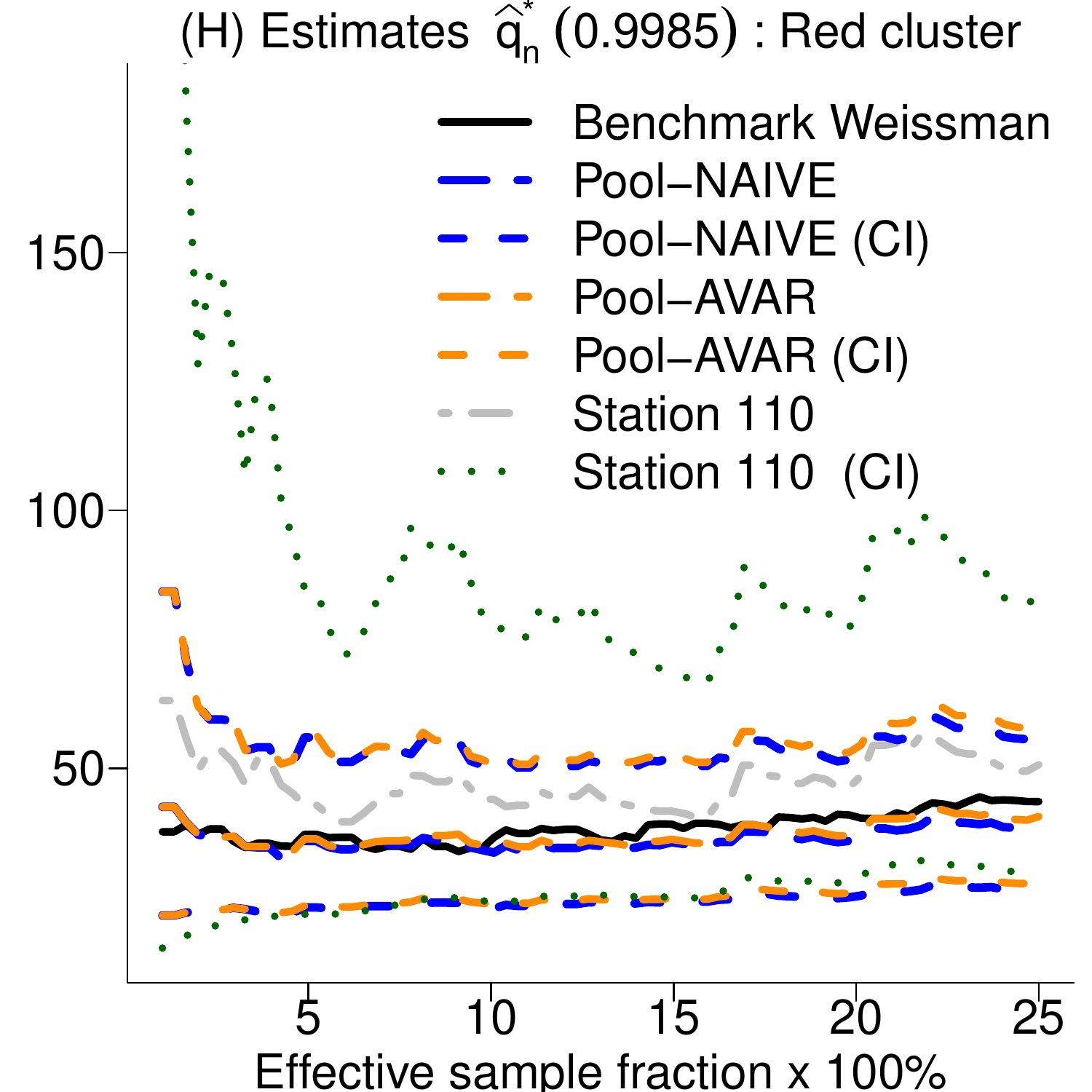}
        \includegraphics[width=0.32\textwidth]{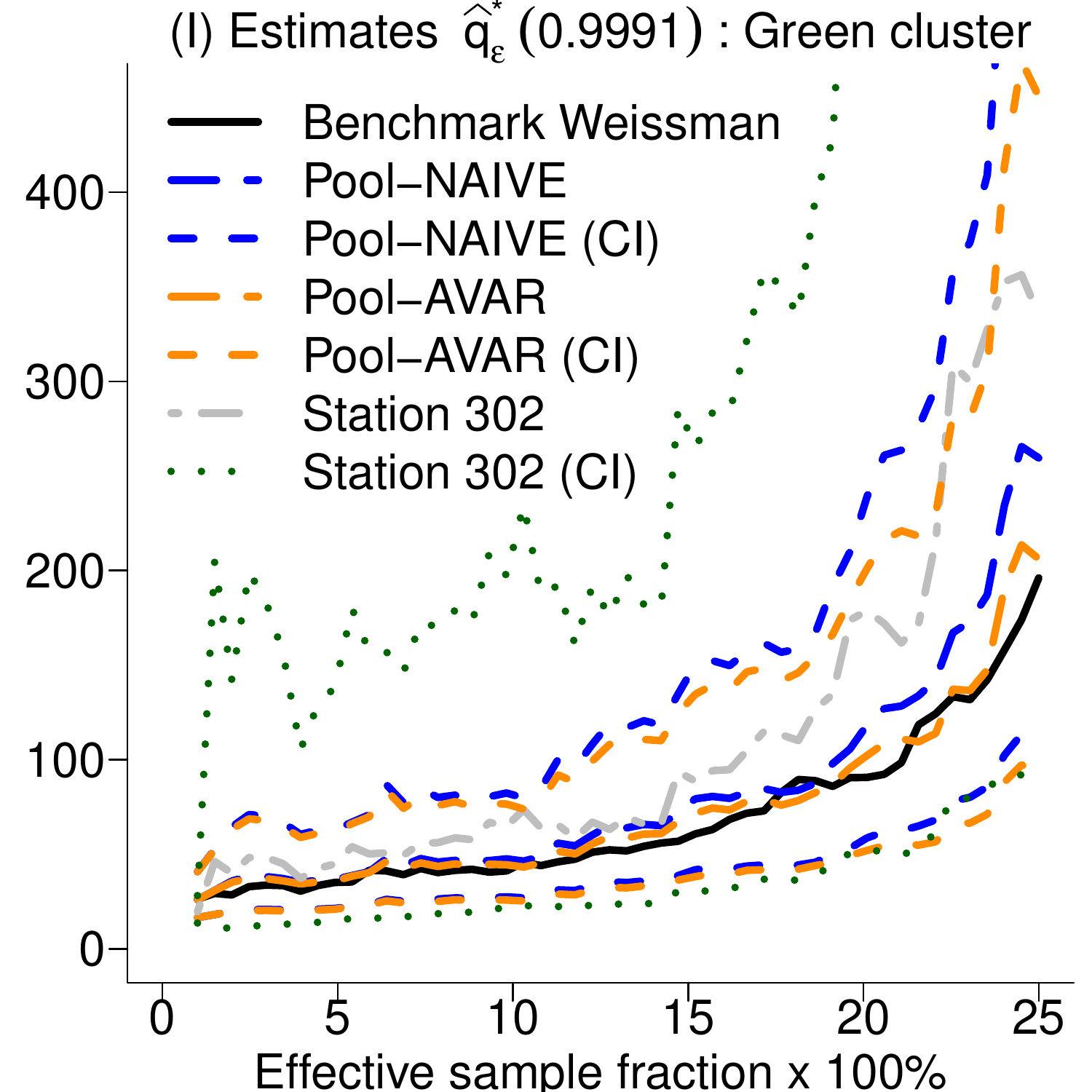} 
\caption{Florida rainfall data. Map of Florida along with its gauge stations and inferential results (notation as in Figure~\ref{fig:carpool:results}).}
\label{fig:florida:Map_test}
\end{figure}

We explore regional inference on tail index and extreme quantiles of monthly rainfall in the state of Florida. An accurate assessment of these tail quantities is crucial for effective flood protection at minimal ecological damage and economic cost. Rainfall measurements are collected daily by the Florida Automated Weather Network at 49 gauge stations, over different periods between December 1997 and May 2021\footnote{See \url{https://fawn.ifas.ufl.edu/data/fawnpub/}.}.
We focus on the eight stations indicated with pin markers in the map in the top panel of Figure~\ref{fig:florida:Map_test}, whose aggregated monthly rainfall exhibit heavy-tailed distributions. The upper tail heaviness for each sample was ascertained in an exploratory analysis using moment and generalized Hill estimators (see {\it e.g.}~\cite{beigoesegteu2004}). 
Individual sample sizes $n_j$ are rather short however, ranging from 172 to 281. The standard extreme value practice of individual extreme value inference will thus be subject to large uncertainty because of the limited amount of data 
at each site. By contrast, our optimal weighted pooling approach allows to reduce the uncertainty by borrowing tail information from homogeneous stations.
%

Our exploratory analysis shows first that the 
eight monthly time series are all stationary according to the classical KPSS and ADF tests. 
A first distinctive property is that the data for each station in the cluster of red pinned stations near the northern border of Florida are, in contrast to those in the green cluster close to the east coast, not autocorrelated according to the Ljung-Box test. The five time series $(X_t)$ in the green cluster of stations stretching along the east coast can be nicely fitted by simple seasonal ARMA models 
$\Phi_{12}(B^{12}) \Phi(B) X_t = c+ \Theta_{12}(B^{12}) \Theta(B) \varepsilon_t$,
where $\varepsilon_t$ stands for a white noise process and $c$ denotes a constant, with $\Phi(B)$, $\Theta(B)$, $\Phi_{12}(B^{12})$ and $\Theta_{12}(B^{12})$ being respectively polynomials of degree $p,q\in\{0,1,2\}$ in the lag operator $B$ and $p_{12},q_{12}\in\{1,2\}$ in $B^{12}$.
Following our theory in Section~\ref{sec:filtering}, the obtained residuals $(\widehat\varepsilon_{t,j}^{(n_j)}, 1\le t\le n_j)$, for each station $j=1,\dots,5$, 
are the basic tool for estimating tail indices. 
A second distinctive property is that the three stations in the red cluster have very similar Hill estimates between $0.32$ and $0.34$, while the 
five stations in the green cluster also have similar Hill estimates between $0.43$ and $0.51$, that are rather different from those in the red cluster. 

This motivates testing for equal tail indices in the separate red and green clusters of stations. First, that the eight stations do not have the same tail index is confirmed by the tail homogeneity test in Figure~\ref{fig:florida:Map_test}(A), where the plot of p-values from the test statistic~$\Lambda_n$ becomes clearly very stable 
below the three significance levels $0.10,\, 0.05$ and~$0.01$. By contrast, we can comfortably conclude the equality of tail indices in both red and green clusters from the tail homogeneity tests in Figure~\ref{fig:florida:Map_test}(B) and~(C). This justifies the in-group tail homogeneity between the stations in each cluster besides their geographical proximity. Even more strongly, the tail homoskedasticity test implemented in Figure~\ref{fig:florida:Map_test}(D) and~(E), for the two extreme quantile levels $1-p\approx 1-1/\max_j n_j$ and $1-p\approx 1-1/n$, allows to accept the assumption of extreme quantile equivalence across stations in each cluster. Therefore, the Hill and Weissman estimators of the common tail index $\gamma$ and extreme quantiles $q(1-p)$, respectively, could be directly calculated from the combined data in each cluster. However, it should be clear that the key question of inference based on these ideal estimators remains 
open in this particular application. Indeed, combining subsamples in each cluster of stations results in a single sample of asymptotically dependent (rainfall) data for which the asymptotic theory of both Hill and Weissman estimators is still unavailable in the extreme value literature. 
Our regional pooled 
estimators come with a satisfactory solution, reducing substantially the huge uncertainty inherent to local inference at each site, as can be seen from Figure~\ref{fig:florida:Map_test}(F)-(I). 
For the red cluster, Figure~\ref{fig:florida:Map_test}(F) shows that both naive and variance-optimal pooled estimators of $\gamma$ are very close to the benchmark Hill estimator, while the asymptotic $95\%$ variance-optimal confidence
intervals are quite stable and remarkably narrower relative to the Hill-based confidence intervals obtained individually from each subsample. We arrive at the same conclusion for the green cluster in Figure~\ref{fig:florida:Map_test}(G), where both 
pooling-type confidence intervals appear to be much tighter 
than the individual Hill-based confidence interval obtained from the largest subsample. Likewise, when estimating the extreme quantile of order $1-p\approx 1-1/n$, the individual Weissman-based confidence interval obtained from the largest subsample in Figure~\ref{fig:florida:Map_test}(H) and~(I), for raw data in the red cluster and for residuals in the green cluster, tends to be unstable and twice as wide as our 
pooling confidence intervals, owing to the reduction of uncertainty in the latter. It is also worth noticing that the variance-optimal pooled estimator is closer to the benchmark Weissman estimator than the naive pooled competitor.

\bibliographystyle{chicago}
\bibliography{reflist}

\end{document}